\documentclass[12pt,reqno]{amsart}
\usepackage{amsthm,amsfonts,amssymb,euscript}

\newcommand{\R}{\ensuremath{\mathbb{R}}}
\newcommand{\Z}{\ensuremath{\mathbb{Z}}}
\newcommand{\C}{\ensuremath{\mathbb{C}}}
\newcommand{\Ss}{\ensuremath{\mathcal{S}}}
\newcommand{\e}{\ensuremath{\mathbf{e}}}
\newcommand{\jma}{\ensuremath{\max(j_1,j_2)}}
\newcommand{\jmi}{\ensuremath{\min(j_1,j_2)}}

\newcommand{\FF}{\ensuremath{\mathcal{F}}}

\numberwithin{equation}{section}
\begin{document}

\title[Global well-posedness of Schr\"{o}dinger maps]
{Low-regularity Schr\"{o}dinger maps, II: global well-posedness in dimensions $d\geq 3$}
\author{Alexandru D. Ionescu}
\address{University of Wisconsin--Madison}
\email{ionescu@math.wisc.edu}
\author{Carlos E. Kenig}
\address{University of Chicago}
\email{cek@math.uchicago.edu}
\thanks{The first author was
supported in part by an NSF grant and a Packard fellowship.
The second author was supported in part by an NSF grant.}
\begin{abstract}
In dimensions $d\geq 3$, we prove that the Schr\"{o}dinger map initial-value problem
\begin{equation*}
\begin{cases}
&\partial_ts=s\times\Delta_x s\,\text{ on }\,\mathbb{R}^d\times\mathbb{R};\\
&s(0)=s_0
\end{cases}
\end{equation*}
is globally well-posed for small data $s_0$ in the critical Besov spaces $\dot{B}_Q^{d/2}(\mathbb{R}^d;\mathbb{S}^2)$, $Q\in\mathbb{S}^2$.
\end{abstract}
\maketitle\tableofcontents

\section{Introduction}\label{intro}

We consider the Schr\"{o}dinger map initial-value problem
\begin{equation}\label{Sch1}
\begin{cases}
&\partial_ts=s\times\Delta_x s\,\text{ on }\,\mathbb{R}^d\times\mathbb{R};\\
&s(0)=s_0,
\end{cases}
\end{equation}
where $d\geq 3$ and
$s:\mathbb{R}^d\times\mathbb{R}\to\mathbb{S}^2\hookrightarrow\mathbb{R}^3$
is a continuous function. The Schr\"{o}dinger map equation has a
rich geometric structure and arises naturally in a number of
different ways; we refer the reader to \cite{NaStUh} or
\cite{KePoStTo} for details. In this paper, which is a
continuation of our earlier work \cite{IoKe2}, we prove a global
well-posedness result for the initial-value problem \eqref{Sch1}
for small data in the critical Besov spaces
$\dot{B}_Q^{d/2}(\mathbb{R}^d;\mathbb{S}^2)$ defined below.

Let $C(\mathbb{R}^d)=\{f:\mathbb{R}^d\to\C:f\text{ is continuous and bounded}\}$. For $\sigma\geq d/2$ we define the Besov-type spaces\footnote{For $\sigma>d/2$ one may replace the space $\dot{B}^\sigma$ with $\dot{B}^{d/2}\cap\dot{H}^\sigma$ throughout the paper (only minor changes would be needed in section \ref{proof}), where $\dot{H}^\sigma$ is the usual homogeneous Sobolev space. We use the Besov spaces $\dot{B}^\sigma$ to measure higher smoothness of functions mostly for simplicity of notation.}
\begin{equation*}
\begin{split}
\dot{B}^{\sigma}=\dot{B}^{\sigma}(\R^d)=&\{\phi\in
C(\R^d):\phi=
\lim_{N\to\infty}\sum_{k=-N}^{N}\mathcal{F}^{-1}_{(d)}[\mathcal{F}_{(d)}\phi\cdot
\eta_k^{(d)}(\xi)]\,\text{ and }\\
&\|\phi\|_{\dot{B}^{\sigma}}=\sum_{k\in\Z}\max(2^{dk/2},2^{\sigma
k})\cdot\|\FF_{(d)}(\phi)(\xi)\cdot \eta_k^{(d)}(\xi
)\|_{L^2}<\infty\},
\end{split}
\end{equation*}
where $\FF_{(d)}$ and $\mathcal{F}^{-1}_{(d)}$ denote the Fourier
transform and the inverse Fourier transform on $\Ss'(\R^d)$, and
$\{\eta_k^{(d)}\}_{k\in\Z}$ is a smooth partition of $1$ with
$\eta_k^{(d)}$ supported in the set
$\{\xi\in\mathbb{R}^d:|\xi|\in[2^{k-1},2^{k+1}]\}$ (see section
\ref{notation} for precise definitions). Let
\begin{equation*}
\dot{B}^\infty=\dot{B}^{\infty}(\mathbb{R}^d)=\bigcap_{\sigma\geq
d/2}\dot{B}^{\sigma}(\mathbb{R}^d)\text{ with the induced metric.}
\end{equation*}

For $\sigma\in[d/2,\infty]$ and $Q=(Q_1,Q_2,Q_3)\in\mathbb{S}^2$
we define the complete metric spaces
\begin{equation}\label{Sch3}
\dot{B}_Q^{\sigma}(\mathbb{R}^d;\mathbb{S}^2)=\{f:\mathbb{R}^d\to\mathbb{R}^3:|f(x)|\equiv
1\text{ and }f_l-Q_l\in \dot{B}^{\sigma}\text{ for }l=1,2,3\},
\end{equation}
with the induced distance
\begin{equation}\label{Ban3}
d^{\sigma}(f,g)=\sum_{l=1}^3\|f_l-g_l\|_{\dot{B}^{\sigma}}.
\end{equation}

For $Q\in \mathbb{S}^2$ let $f_Q(x)\equiv Q$, $f_Q\in
H^{\infty}_Q(\mathbb{R}^d;\mathbb{S}^2)$. For any metric space
$X$, $x\in X$, and $r>0$ let $B_X(x,r)$ denote the open ball
$\{y\in X:d(x,y)<r\}$. Let $\mathbb{Z}_+=\{0,1,\ldots\}$. Our main
theorem concerns global well-posedness of the initial-value
problem \eqref{Sch1} for small data $s_0\in
\dot{B}_Q^{d/2}(\mathbb{R}^d;\mathbb{S}^2)$, $Q\in\mathbb{S}^2$.

\newtheorem{Main1}{Theorem}[section]
\begin{Main1}\label{Main1}
(a) Assume $d\geq 3$ and $Q\in\mathbb{S}^2$. Then there are
numbers $\epsilon_0\leq\overline{\epsilon}_0\in(0,\infty)$ with
the property that for any $s_0\in
\dot{B}^{\infty}_Q(\mathbb{R}^d;\mathbb{S}^2)\cap
B_{\dot{B}^{d/2}_Q(\mathbb{R}^d;\mathbb{S}^2)}(f_Q,\epsilon_0)$
there is a unique solution
\begin{equation*}
s=S^\infty_Q(s_0)\in
C(\R:\dot{B}^{\infty}_Q(\mathbb{R}^d;\mathbb{S}^2)\cap
B_{\dot{B}^{d/2}_Q(\mathbb{R}^d;\mathbb{S}^2)}(f_Q,\overline{\epsilon}_0))
\end{equation*}
of the initial-value problem \eqref{Sch1}.

(b) In addition, we have the Lipschitz bound
\begin{equation}\label{ak1}
\sup_{t\in\R}d^{d/2}(S^{\infty}_Q(s_0)(t),S^{\infty}_Q(s'_0)(t))\leq
C\cdot d^{d/2}(s_0,s'_0)
\end{equation}
for  any $s_0,s'_0\in
\dot{B}^{\infty}_Q(\mathbb{R}^d;\mathbb{S}^2)\cap
B_{\dot{B}^{d/2}_Q(\mathbb{R}^d;\mathbb{S}^2)}(f_Q,\epsilon_0)$.
Thus the mapping $s_0\to S^\infty_Q(s_0)$ extends uniquely to a
Lipschitz mapping
\begin{equation*}
S^{d/2}:B_{\dot{B}^{d/2}_Q(\mathbb{R}^d;\mathbb{S}^2)}(f_Q,\epsilon_0)\to
C(\R:B_{\dot{B}^{d/2}_Q(\mathbb{R}^d;\mathbb{S}^2)}(f_Q,\overline{\epsilon}_0)).
\end{equation*}
\end{Main1}

Theorem \ref{Main1} appears to be the first low-regularity global well-posedness result for the Schr\"o\-din\-ger map initial-value problem.  Its direct analogue in the setting of wave maps is the work of Tataru \cite{Tat1} (see also \cite{KlMa}, \cite{KlSe}, \cite{Tat2}, \cite{Ta1}, \cite{Ta2}, \cite{KlRo}, \cite{ShSt}, and \cite{Tat3} for other local and global well-posedness theorems for wave maps).

The initial-value problem \eqref{Sch1} has been studied extensively (also in the case in which the sphere $\mathbb{S}^2$ is replaced by more general targets). It is known that sufficiently smooth solutions exist locally in time, even for large data (see, for example, \cite{SuSuBa}, \cite{ChShUh}, \cite{DiWa2}, \cite{Ga}, \cite{KePoStTo} and the references therein). Such theorems for (local in time) smooth solutions are proved using variants of the energy method. For low-regularity data, the energy method cannot be applied, and the initial-value problem \eqref{Sch1} has been studied indirectly using the ``modified Schr\"{o}dinger map equation'' (see, for example, \cite{NaStUh}, \cite{NaStUh2}, \cite{KeNa}, and \cite{KaKo}). While existence and uniqueness theorems for this modified Schr\"{o}dinger map equation in certain low-regularity spaces are known (at least in dimension $d=2$), it is not clear whether such theorems can be transfered to the original Schr\"{o}dinger map initial-value problem (see, however, \cite{NaShVeZe}). 

In \cite{IoKe2}, the authors proved local well-posedness of the initial-value problem \eqref{Sch1} for small data in the natural Sobolev spaces $H^\sigma_Q(\mathbb{R}^d;\mathbb{S}^2)$, $\sigma>(d+1)/2$. This was achieved by reducing the initial-value problem \eqref{Sch1} to the nonlinear Schr\"{o}dinger equation \eqref{Sch5} below\footnote{Using the stereographic projection, such a reduction is possible due to the fact that the solutions take values only in a small part of the sphere; the models \eqref{Sch1} and \eqref{Sch5} are certainly not equivalent without such a smallness assumption.}, and by analyzing the resulting equation using a direct perturbative argument. We follow the same approach in this paper. At about the same time and independently, Bejenaru \cite{Be} proved local well-posedness of the initial-value  problem  \eqref{Sch5} for small data in the Sobolev spaces $H^\sigma$, for $\sigma$ in the full subcritical range $\sigma>d/2$. The resolution spaces used by Bejenaru \cite{Be} appear to be very different from the spaces used by us in \cite{IoKe2} and in this paper. To reach the full subcritical range, Bejenaru noticed, apparently for the first time in the setting of Schr\"{o}dinger maps, that the gradient part of the nonlinearity  in \eqref{Sch5} has a certain null structure (similar to the null structure of the wave maps). We exploit this null structure through the identity  \eqref{null}.  

Theorem \ref{Main1} can be restated using the stereographic projection.
By rotation invariance, we may assume
\begin{equation}\label{al1}
Q=(0,0,1).
\end{equation}
Assume $\epsilon>0$ is small enough. For $f=(f_1,f_2,f_3)\in
B_{\dot{B}^{d/2}_Q(\mathbb{R}^d;\mathbb{S}^2)}(f_Q,\epsilon)$ let
\begin{equation*}
g=L(f)=\frac{f_1+if_2}{1+f_3}.
\end{equation*}
Clearly, $L(f):\mathbb{R}^d\to\mathbb{C}$ is continuous and takes
values in a small neighborhood of $0$. For $g\in
B_{\dot{B}^{d/2}}(0,\epsilon)$ we define
\begin{equation*}
f=(f_1,f_2,f_3)=\widetilde{L}(g)=\Big(\frac{g+\overline{g}}
{1+g\overline{g}},\frac{(-i)(g-\overline{g})}
{1+g\overline{g}},\frac{1-g\overline{g}}{1+g\overline{g}}\Big).
\end{equation*}
Clearly, $\widetilde{L}(g):\mathbb{R}^d\to\mathbb{S}^2$ is
continuous and takes values in a small neighborhood of $Q$. A
direct computation shows that $u:\mathbb{R}^d\to\{z\in\C:|z|\leq
1\}$ is a smooth solution of the equation
\begin{equation*}
(i\partial_t+\Delta_x)u=\frac{2\overline{u}}{1+u\overline{u}}\sum_{j=1}^d(\partial_{x_j}u)^2\text{
on }\mathbb{R}^d\times\R,
\end{equation*}
if and only if the function
$s:\mathbb{R}^d\to\mathbb{S}^2\cap\{(x_1,x_2,x_3)\in\R^3:x_3\geq
0\}$, $s(t)=\widetilde{L}(u(t))$, is a smooth solution of the
Schr\"{o}dinger map equation
\begin{equation*}
\partial_ts=s\times\Delta_x s\,\text{ on }\,\mathbb{R}^d\times\R.
\end{equation*}
Since $\dot{B}^\sigma$, $\sigma\in[d/2,\infty)$ are Banach
algebras, in the sense that $$||uv||_{\dot{B}^\sigma}\leq
C_\sigma(||u||_{\dot{B}^\sigma}||v||_{\dot{B}^{d/2}}
+||u||_{\dot{B}^{d/2}}||v||_{\dot{B}^\sigma})$$ for any
$\sigma\in[d/2,\infty)$ and $u,v\in\dot{B}^\sigma$, for Theorem
\ref{Main1} it suffices to prove Theorem \ref{Main2} below.

\newtheorem{Main2}[Main1]{Theorem}
\begin{Main2}\label{Main2}
(a) Assume $d\geq 3$. Then there are numbers $\epsilon_1\leq
\overline{\epsilon}_1\in(0,\infty)$ with the property that for any
$\phi \in \dot{B}^\infty \cap B_{\dot{B}^{d/2}}(0,\epsilon_1)$
there is a unique solution
\begin{equation*}
u=\widetilde{S}^\infty(\phi)\in C(\R:\dot{B}^\infty\cap
B_{\dot{B}^{d/2}}(0,\overline{\epsilon}_1))
\end{equation*}
of the initial-value problem
\begin{equation}\label{Sch5}
\begin{cases}
&(i\partial_t+\Delta_x)u=2\overline{u}(1+u\overline{u})^{-1}\sum_{j=1}^d(\partial_{x_j}u)^2\text{ on }\mathbb{R}^d\times\R;\\
&u(0)=\phi.
\end{cases}
\end{equation}

(b) In addition, we have the Lipschitz bound
\begin{equation}\label{by1}
\sup_{t\in\R}\|\widetilde{S}^\infty(\phi)(t)-\widetilde{S}^\infty(\phi')(t)\|_{\dot{B}^{d/2}}\leq
C\|\phi-\phi'\|_{\dot{B}^{d/2}},
\end{equation}
for any $\phi,\phi'\in \dot{B}^\infty\cap
B_{\dot{B}^{d/2}}(0,\epsilon_1)$. Thus the mapping $\phi\to
S^\infty(\phi)$ extends uniquely to a Lipschitz mapping
\begin{equation*}
\widetilde{S}^{d/2}:B_{\dot{B}^{d/2}}(0,\epsilon_1)\to
C(\R:B_{\dot{B}^{d/2}}(0,\overline{\epsilon}_1)).
\end{equation*}
\end{Main2}

By scale invariance, it suffices to construct the solution
$u=S^\infty(\phi)$ on the time interval $[-1,1]$ and prove the Lipschitz 
bound \eqref{by1} for $t\in[-1,1]$. The resolution spaces we construct in section \ref{notation} are adapted to this restriction in time. This restriction creates a somewhat artificial distinction between frequencies that are $\leq 1$ and frequencies that are $\geq 1$. The benefit of this  time restriction, however, is that  the denominators in formulas such as \eqref{gu45.1} and  \eqref{pp2} (and in many other places) do not vanish, and all of our  integrals are absolutely convergent (in particular, changes  of order of integration are justified). The direct use of scale-invariant spaces would lead to denominators such as $\tau+|\xi|^2$, and the integrals containing such denominators would  not converge absolutely. 

The rest of the  paper is organized as follows: in section \ref{notation} we define our main (dyadic) resolution spaces and establish some of their basic properties. These spaces are minor modifications of the resolution spaces already used by the authors in \cite{IoKe2} (see also \cite{IoKe} for the one-dimensional analogues of these resolution spaces). In section \ref{proof} we give the main argument that proves Theorem \ref{Main2}; the main ingredients in our perturbative argument are the four nonlinear estimates \eqref{amy1}, \eqref{amy2}, \eqref{amy3}, and  \eqref{amy4}. In the remaining sections we prove these four nonlinear estimates. The key ingredients in these proofs are the scale-invariant $L^{2,\infty}_{\e'}$  (maximal function) estimate in Lemma \ref{Lemmaa1} and the scale-invariant $L^{\infty,2}_{\e'}$ (local smoothing) estimate in Lemma \ref{Lemmas2}. These two estimates have been used before by the authors in \cite{IoKe} and \cite{IoKe2}. The maximal function bound fails (logarithmically) in dimension $d=2$, which  is  the main reason why we need to assume $d\geq 3$. 

We would like to thank Bejenaru for making his preprint \cite{Be} available to us.

\section{Notation and preliminary lemmas}\label{notation}

In this section we summarize most of the notation, define our main normed spaces,\footnote{It is likely that only minor changes would be needed to guarantee that all of our normed spaces are in fact Banach spaces. We do not need this, however, since the limiting argument in the construction of solutions takes place in the Banach space $C(\R:\dot{B}^{d/2})$.}  and prove some of their basic properties. Let $\mathcal{F}$ and $\mathcal{F}^{-1}$ denote the Fourier transform and the inverse Fourier transform operators on $\mathcal{S}'(\mathbb{R}^{d+1})$. For $l=1,\ldots,d$ let $\mathcal{F}_{(l)}$ and $\mathcal{F}_{(l)}^{-1}$ denote the Fourier transform and the inverse Fourier transform operators on $\mathcal{S}'(\mathbb{R}^l)$.

We fix $\eta_0:\mathbb{R}\to[0,1]$ a smooth even function supported in the set $\{\mu\in\mathbb{R}:|\mu|\leq 8/5\}$ and equal to $1$ in the set $\{\mu\in\mathbb{R}:|\mu|\leq 5/4\}$. Then we define $\eta_j:\mathbb{R}\to[0,1]$, $j=1,2,\ldots$,
\begin{equation}\label{gu1}
\eta_j(\mu)=\eta_0(\mu/2^j)-\eta_0(\mu/2^{j-1}),
\end{equation}
and $\eta_k^{(d)}:\mathbb{R}^d\to[0,1]$, $k\in\Z$,
\begin{equation}\label{gu1.1}
\eta_k^{(d)}(\xi )=\eta_0( |\xi|/2^k)-\eta_0( |\xi|/2^{k-1}).
\end{equation}
For $j_1,j_2\in\Z$, we also define $\eta_{[j_1,j_2]}=\sum_{j_1\leq j'\leq j_2}\eta_{j'}$ (with the conventions $\eta_{[j_1,j_2]}\equiv 0$ if $j_1>j_2$ and  $\eta_j\equiv 0$ if $j\leq -1$), $\eta_j^{\pm}(\mu)=\eta_j(\mu)\cdot\mathbf{1}_{[0,\infty)}(\pm\mu)$, $\eta_{[j_1,j_2]}^{\pm}(\mu)=\eta_{[j_1,j_2]}(\mu)\cdot\mathbf{1}_{[0,\infty)}(\pm\mu)$, $\eta_{\leq j}=\eta_{[0,j]}$, $\eta_{\geq j}=1-\eta_{\leq j-1}$.

For $k\in\mathbb{Z}$ let $I_k^{(d)}=\{\xi\in\R^d:|\xi|\in[2^{k-1},2^{k+1}]\}$; for $j\in\Z_+$ let $I_j=\{\mu\in\R:|\mu|\in[2^{j-1},2^{j+1}]\}$ if $j\geq 1$ and $I_j=[-2,2]$ if $j=0$. For $k\in\Z$ and $j\in\Z_+$ let
\begin{equation*}
D_{k,j}=\{(\xi,\tau)\in\mathbb{R}^d\times\mathbb{R}:\xi \in I_k^{(d)}\text{ and }|\tau+|\xi|^2|\in I_j\}\text{ and }D_{k,\leq j}=\bigcup\limits_{0\leq j'\leq j}D_{k,j'}.
\end{equation*}

For $k\in\Z$ we define first the normed spaces
\begin{equation}\label{v1}
\begin{split}
X_k=\{f\in L^2(\mathbb{R}^d\times&\mathbb{R}):f\text{ supported in }I_k^{(d)}\times\mathbb{R}\text { and } \\
&\|f\|_{X_k}=\sum_{j=0}^\infty 2^{j/2}\beta_{k,j}\|\eta_j(\tau+|\xi|^2)\cdot f\|_{L^2}<\infty\},
\end{split}
\end{equation}
where, with $k_+=\max(k,0)$,
\begin{equation}\label{v1.1}
\beta_{k,j}=1+2^{(j-2k_+)/2}.
\end{equation}

The spaces $X_k$ are not sufficient for our estimates, due to various logarithmic divergences. For any vector $\mathbf{e}\in\mathbb{S}^{d-1}$ let $$P_{\mathbf{e}}=\{\xi\in\mathbb{R}^d:\xi\cdot\mathbf{e}=0\}$$ with the induced Euclidean measure. For $p,q\in[1,\infty]$ we define the normed spaces $L^{p,q}_{\mathbf{e}}=L^{p,q}_{\mathbf{e}}(\mathbb{R}^d\times\mathbb{R})$,
\begin{equation}\label{vv1}
\begin{split}
L^{p,q}_{\mathbf{e}}&=\{f\in L^2(\mathbb{R}^d\times \mathbb{R}):\\
&\|f\|_{L^{p,q}_{\mathbf{e}}}=\Big[\int_{\mathbb{R}}\Big[\int_{P_\mathbf{e}\times \mathbb{R}}|f(r\mathbf{e}+v,t)|^q\,dvdt\Big]^{p/q}\,dr\Big]^{1/p}<\infty\}.
\end{split}
\end{equation}
For $k\geq 100$, $j\in \mathbb{Z}_+$ and $k'\in[1,k+1]\cap\Z$ let
\begin{equation*}
D_{k,j}^{\mathbf{e},k'}=\{(\xi,\tau)\in D_{k,j}:\xi\cdot\mathbf{e}\in I_{k'}\cap [0,\infty)\}\text{ and }D_{k,\leq j}^{\mathbf{e},k'}=\bigcup_{0\leq j'\leq j}D_{k,j}^{\mathbf{e},k'}.
\end{equation*}
For $k\geq 100$, $k'\in[1,k+1]\cap\Z$, and $\mathbf{e}\in\mathbb{S}^{d-1}$, we define the normed spaces
\begin{equation}\label{v2}
\begin{split}
Y_k^{\mathbf{e},k'}&=\{f\in L^2(\mathbb{R}^d\times\mathbb{R}):f\text{ supported in }D_{k,\leq 2k+10}^{\mathbf{e},k'}\text { and } \\
&\|f\|_{Y_k^{\mathbf{e},k'}}=2^{-k'/2}\gamma_{k,k'}\cdot \|\mathcal{F}^{-1}[(\tau+|\xi|^2+i)\cdot f]\|_{L^{1,2}_{\mathbf{e}}}<\infty\},
\end{split}
\end{equation}
where
\begin{equation*}
\gamma_{k,k'}=2^{2d(k-k')}.
\end{equation*}
For $k\geq 100$ and $\mathbf{e}\in\mathbb{S}^{d-1}$, we define the normed spaces
\begin{equation}\label{v2.5}
\begin{split}
Y_k^{\mathbf{e}}&=\{f\in L^2(\mathbb{R}^d\times\mathbb{R}):f\text{ supported in }\cup_{k'=1}^{k+1}D_{k,\leq 2k+10}^{\mathbf{e},k'}\text { and } \\
&\|f\|_{Y_k^{\mathbf{e}}}=\sum_{k'=1}^{k+1}\|f\cdot \eta_{k'}^{+}(\xi\cdot\e)\|_{Y_k^{\mathbf{e},k'}}<\infty\}.
\end{split}
\end{equation}
For simplicity of notation, we also define $Y_k^{\mathbf{e}}=\{0\}$ for $k\leq 99$.

We fix $L=L(d)$ large and $\mathbf{e}_1,\ldots,\mathbf{e}_L\in\mathbb{S}^{d-1}$, $\mathbf{e}_l\neq \mathbf{e}_{l'}$ if $l\neq l'$, such that
\begin{equation}\label{vm4}
\text{ for any }\mathbf{e}\in\mathbb{S}^{d-1}\text{ there is }l\in\{1,\ldots,L\}\text{ such that }|\mathbf{e}-\mathbf{e}_l|\leq 2^{-100}.
\end{equation}
We assume in addition that if $\mathbf{e}\in\{\mathbf{e}_1,\ldots,\mathbf{e}_L\}$ then $-\mathbf{e}\in\{\mathbf{e}_1,\ldots,\mathbf{e}_L\}$. For $k\in\mathbb{Z}$ we define the normed spaces
\begin{equation}\label{v3'}
Z_k=X_k+Y_{k}^{\mathbf{e}_1}+\ldots+Y_k^{\mathbf{e}_L}.
\end{equation}
The spaces $Z_k$ are our main normed spaces.

We prove now several estimates. In view of the definitions, if $m\in L^\infty(\mathbb{R}^d)$, $\mathcal{F}_{(d)}^{-1}(m)\in L^1(\mathbb{R}^d)$, and $f\in Z_k$,  then $m(\xi)\cdot f\in Z_k$ and
\begin{equation}\label{mi2}
||m(\xi)\cdot f||_{Z_k}\leq C||\mathcal{F}_{(d)}^{-1}(m)||_{L^1(\mathbb{R}^d)}\cdot ||f||_{Z_k}.
\end{equation}
We show first that the spaces $Z_k$ are logarithmic modifications of the spaces $X_k$.

\newtheorem{Lemmas1}{Lemma}[section]
\begin{Lemmas1}\label{Lemmas1}
If $k\in\Z$, $j\in\mathbb{Z}_+$ and $f\in Z_k$ then
\begin{equation}\label{vvs1}
\|f\cdot \eta_j(\tau+|\xi|^2)\|_{X_k}\leq C\|f\|_{Z_k}.
\end{equation}
\end{Lemmas1}

\begin{proof}[Proof of Lemma \ref{Lemmas1}] Clearly, we may assume $k\geq 100$ and $f=f^{\mathbf{e},k'}\in Y_k^{\mathbf{e}}$, for some $\mathbf{e}\in\{\mathbf{e}_1,\ldots,\mathbf{e}_L\}$ and $k'\in[1,k+1]\cap\Z$. Let
\begin{equation}\label{gu45}
h(x,t)=2^{-k'/2}\mathcal{F}^{-1}[(\tau+|\xi|^2+i)\cdot f^{\mathbf{e},k'}](x,t).
\end{equation}
Thus
\begin{equation}\label{gu45.1}
f^{\mathbf{e},k'}(\xi,\tau)=\mathbf{1}_{D_{k,\leq 2k+10}^{\mathbf{e},k'}}(\xi,\tau)\cdot\frac{2^{k'/2}}{\tau+|\xi|^2+i}\mathcal{F}(h)(\xi,\tau).
\end{equation}
In view of the definitions, for \eqref{vvs1} it suffices to prove the stronger bound
\begin{equation}\label{gu21}
2^{k'/2}2^{-j/2}\|\mathbf{1}_{D_{k,j}^{\e,k'}}(\xi,\tau)\cdot\mathcal{F}(h)\|_{L^2_{\xi,\tau}}\leq C\|h\|_{L^{1,2}_{\mathbf{e}}}
\end{equation}
for any $j\leq 2k+11$. We write $\xi=\xi_1\mathbf{e}+\xi'$, $x=x_1\mathbf{e}+x'$, $x_1,\xi_1\in\mathbb{R}$, $x',\xi'\in P_{\mathbf{e}}$. Let
\begin{equation*}
h'(x_1,\xi',\tau)=\int_{P_{\mathbf{e}}\times\mathbb{R}}h(x_1\mathbf{e}+x',t)e^{-i(x'\cdot \xi'+t\tau)}\,dx'dt.
\end{equation*}
By Plancherel theorem,
\begin{equation*}
\|h\|_{L^{1,2}_{\mathbf{e}}}=C\|h'\|_{L^1_{x_1}L^2_{\xi',\tau}}.
\end{equation*}
Thus, for \eqref{gu21}, it suffices to prove that
\begin{equation}\label{pp1}
\begin{split}
2^{(k'-j)/2}\Big|\Big|\mathbf{1}_{D_{k,j}^{\e,k'}}(\xi,\tau)\cdot&\int_\mathbb{R}h'(x_1,\xi',\tau)e^{-ix_1\xi_1}\,dx_1\Big|\Big|_{L^2_{\xi_1,\xi',\tau}}\\
&\leq C\min(1,2^{k'-j/2})\cdot\|h'\|_{L^1_{x_1}L^2_{\xi',\tau}}.
\end{split}
\end{equation}
This follows easily since for any $(\xi',\tau)\in P_{\mathbf{e}}\times\mathbb{R}$ the measure of the set $\{\xi_1:|\xi_1|\approx 2^{k'},\,|\tau+\xi_1^2+|\xi'|^2|\leq 2^{j+1}\}$ is bounded by $C\min(2^{k'},2^{j-k'})$.
\end{proof}

The implicit gain of $\gamma_{k,k'}$ in the bound \eqref{pp1} shows that $Y_{k}^{\e,k'}\hookrightarrow X_k$ if $k'\leq 9k/10$. Let $T_k=[9k/10,k+1]\cap\Z$. In view of the definitions, if $f\in Z_k$ then we can write
\begin{equation}\label{mi1}
\begin{cases}
\negmedspace\negmedspace\negmedspace\negmedspace\negmedspace&f=\sum\limits_{j\in\mathbb{Z}_+}g_j+\sum\limits_{k'\in T_k}\sum\limits_{l=1}^Lf^{\mathbf{e}_l,k'},\,g_j\text{ supported in }D_{k,j},\,f^{\mathbf{e}_l,k'}\text{ supported in }D_{k,\leq 2k+10}^{\mathbf{e}_l,k'};\\
\negmedspace\negmedspace\negmedspace\negmedspace\negmedspace&\sum\limits_{j\in\mathbb{Z}_+}2^{j/2}\beta_{k,j}||g_j||_{L^2}+\sum\limits_{k'\in T_k}\sum\limits_{l=1}^L||f^{\mathbf{e}_l,k'}||_{Y_k^{\mathbf{e}_l,k'}}\leq 2\|f\|_{Z_k}.
\end{cases}
\end{equation}
This is our main atomic decomposition of functions in $Z_k$.

In addition, the bound \eqref{pp1} shows that if $k\in\Z$ and
\begin{equation*}
f\text{ is supported in }I_k^{(d)}\times\R\cap
\{(\xi,\tau):\xi\cdot\e\geq 2^{k-40}\}\text{ for some }\e\in\mathbb{S}^{d-1},
\end{equation*}
then, for any $j\geq 0$,
\begin{equation}\label{xx1}
||f\cdot\eta_{j}(\tau+|\xi|^2)||_{Z_k}\leq C 2^{-k/2}||\mathcal{F}^{-1}
[(\tau+|\xi|^2+i)\cdot f]||_{L^{1,2}_\e}.
\end{equation}

We prove now $L^\infty_tL^2_x$ and $L^\infty_{x,t}$ estimates.

\newtheorem{Lemmaa2}[Lemmas1]{Lemma}
\begin{Lemmaa2}\label{Lemmaa2}
If $k\in\Z$, $t\in\mathbb{R}$, and $f\in Z_k$ then
\begin{equation}\label{lb4}
\sup_{t\in\mathbb{R}}\|\mathcal{F}^{-1}(f)(.,t)\|_{L^2_x}
\leq C\|f\|_{Z_k}.
\end{equation}
Thus
\begin{equation}\label{lb44}
\|\mathcal{F}^{-1}(f)\|_{L^\infty_{x,t}}
\leq C2^{dk/2}\|f\|_{Z_k}.
\end{equation}
\end{Lemmaa2}

\begin{proof}[Proof of Lemma \ref{Lemmaa2}] By Plancherel theorem it
suffices to prove that
\begin{equation}\label{si1}
\Big|\Big|\int_{\mathbb{R}}f(\xi,\tau)e^{it\tau}\,d\tau\Big|\Big|_{L^2_\xi}\leq C\|f\|_{Z_k}.
\end{equation}
We use the representation \eqref{mi1}. Assume first that $f=g_j$. Then
\begin{equation}\label{si1.1}
\begin{split}
\Big|\Big|\int_\mathbb{R}g_j(\xi,\tau)e^{it\tau}\,d\tau\Big|\Big|_{L^2_\xi}\leq C||g_j(\xi,\tau)||_{L^2_\xi L^1_\tau}\leq C2^{j/2}||g_j||_{L^2_{\xi,\tau}},
\end{split}
\end{equation}
which proves \eqref{si1} in this case. This inequality also shows that
\begin{equation}\label{si1.2}
\|\mathcal{F}^{-1}(g_j)\|_{L^\infty}\leq C2^{dk/2}2^{j/2}\|g_j\|_{L^2}.
\end{equation}

Assume now that $k\geq 100$ and $f=f^{\mathbf{e},k'}\in Y_k^{\mathbf{e},k'}$, $\mathbf{e}\in\{\mathbf{e}_1,\ldots,\mathbf{e}_L\}$, $k'\in T_k$. We have to prove that
\begin{equation}\label{si1'}
\Big|\Big|\int_{\mathbb{R}}f^{\mathbf{e},k'}(\xi,\tau)e^{it\tau}\,d\tau\Big|\Big|_{L^2_\xi}\leq C\|f^{\mathbf{e},k'}\|_{Y_k^{\mathbf{e},k'}}.
\end{equation}
We define the function $h$ as in \eqref{gu45}. In view of \eqref{pp1},
\begin{equation}\label{vc1}
\|\eta_{\geq 2k'-49}(\tau+|\xi|^2)\cdot f^{\mathbf{e},k'}\|_{X_k}\leq C\|f^{\mathbf{e},k'}\|_{Y^{\e,k'}_k}.
\end{equation}
Since the bound \eqref{si1'} was already proved for functions in $X_k$ (see \eqref{si1.1}), for \eqref{si1'} it suffices to prove the stronger bound
\begin{equation}\label{vc2}
\Big|\Big|\int_{\mathbb{R}}f^{\mathbf{e},k'}(\xi,\tau)\cdot \eta_{\leq 2k'-50}(\tau+|\xi|^2)\cdot e^{it\tau}\,d\tau\Big|\Big|_{L^2_\xi}\leq C\|f^{\mathbf{e},k'}\|_{Y_k^{\mathbf{e},k'}}.
\end{equation}
We use the formula \eqref{gu45.1}, and write $\xi=\xi_1\e+\xi'$, $\xi_1\in\R$, $\xi'\in P_\e$. For \eqref{vc2} it suffices to prove that
\begin{equation}\label{si2}
2^{k'/2}\Big|\Big|\eta_{[k'-1,k'+1]}^+(\xi_1)\int_{\mathbb{R}}\frac{\eta_{\leq 2k'-50}(\tau+|\xi|^2)}{\tau+|\xi|^2+i}\cdot \mathcal{F}(h)(\xi,\tau)\cdot e^{it\tau}\,d\tau\Big|\Big|_{L^2_\xi}\leq C||h||_{L^{1,2}_\e},
\end{equation}
for any $t\in\mathbb{R}$.
As in Lemma \ref{Lemmas1}, for \eqref{si2} it suffices to prove that
\begin{equation}\label{vc6}
2^{k'/2}\Big|\Big|\eta_{[k'-1,k'+1]}^+(\xi_1)\int_{\mathbb{R}}\frac{\eta_{\leq 2k'-50}(\tau+|\xi|^2)}{\tau+|\xi|^2+i}\cdot h'(\xi',\tau)\cdot e^{it\tau}\,d\tau\Big|\Big|_{L^2_\xi}\leq C||h'||_{L^2_{\xi',\tau}},
\end{equation}
for any $t\in\R$ and $h'\in L^2(P_\e\times\R)$. We may assume $t=0$ and let
\begin{equation*}
h''(\xi',\mu)=\int_{\mathbb{R}}\frac{\eta_{\leq 2k'-50}(\tau+\mu)}{\tau+\mu+i}h'(\xi',\tau)\,d\tau.
\end{equation*}
In view of the boundedness of the Hilbert transform on $L^2(\mathbb{R})$, $||h''||_{L^2_{\xi',\mu}}\leq C||h'||_{L^2_{\xi',\tau}}$. Thus, for \eqref{vc6}, it suffices to prove that
\begin{equation*}
2^{k'/2}||\eta_{[k'-1,k'+1]}^+(\xi_1)\cdot h''(\xi',\xi_1^2+|\xi'|^2)||_{L^2_{\xi',\xi_1}}\leq C||h''||_{L^2_{\xi',\mu}},
\end{equation*}
which follows easily by changes of variables.
\end{proof}

We consider now the action of multipliers  of the form $m_{\leq j}(\tau+|\xi|^2)$.

\newtheorem{Lemmag1}[Lemmas1]{Lemma}
\begin{Lemmag1}\label{Lemmag1}
Assume $m:\mathbb{R}\to\C$ is a smooth function supported in the interval $[-2,2]$ and let $m_{ \leq j}(\mu)=m(\mu/2^j)$, $j\in\Z_+$. If $k\geq 100$, $k'\in T_k$, $j\in[0,2k'-80]\cap\Z$, $\e\in\mathbb{S}^{d-1}$, and $f$ is supported in $I_k^{(d)}\times\R$ then
\begin{equation}\label{hv1}
||\mathcal{F}^{-1}[m_{\leq  j}(\tau+|\xi|^2)\cdot f\cdot \eta_{k'}^+(\xi\cdot\e)]||_{L^{1,2}_\e}\leq C||\mathcal{F}^{-1}(f)||_{L^{1,2}_\e}.
\end{equation}
Thus, if $k\in\Z$, $j\in\Z_+$, and $f\in Z_k$ then
\begin{equation}\label{bt70}
||\eta_{\leq j}(\tau+|\xi|^2)\cdot f||_{Z_k}\leq C||f||_{Z_k}.
\end{equation}
\end{Lemmag1}

\begin{proof} We write as before $\xi=\xi_1\mathbf{e}+\xi'$, $x=x_1\mathbf{e}+x'$, $x_1,\xi_1\in\mathbb{R}$, $x',\xi'\in P_{\mathbf{e}}$. Using  Plancherel theorem, it suffices to prove that
\begin{equation}\label{hv2}
\Big|\Big|\eta_{\leq 2k+10}(\tau+|\xi'|^2)\int_{\mathbb{R}}e^{ix_1\xi_1}m_{\leq  j}(\tau+|\xi'|^2+\xi_1^2)\cdot \eta_{k'}^+(\xi_1)\,d\xi_1\Big|\Big|_{L^1_{x_1}L^\infty_{\xi',\tau}}\leq C.
\end{equation}
In view of the restriction $j\leq 2k'-80$, we  may assume that the supremum in $(\xi',\tau)$ in \eqref{hv2} is taken over the set $\{(\xi',\tau):-\tau-|\xi'|^2\in[2^{2k'-70},2^{2k'+10}]\}$. Let $M=M(\xi',\tau)=(-\tau-|\xi'|^2)^{1/2}$, $M\approx  2^{k'}$. By integration by parts
\begin{equation*}
\Big|\int_{\mathbb{R}}e^{ix_1\xi_1}m_{\leq  j}(\xi_1^2-M^2)\cdot \eta_{k'}^{+}(\xi_1)\,d\xi_1\Big|\leq C\frac{2^{j-k'}}{1+(2^{j-k'}x_1)^2},
\end{equation*}
if $M\approx 2^{k'}$, which gives  \eqref{hv2}.

The inequality \eqref{bt70} follows from \eqref{hv1} and \eqref{pp1}.
\end{proof}

We conclude this section with a representation formula for functions in $Y_k^{\e,k'}$.

\newtheorem{Lemmar1}[Lemmas1]{Lemma}
\begin{Lemmar1}\label{Lemmar1}
If $k\geq 100$, $k'\in T_k$, $\mathbf{e}\in\{\e_1,\ldots,\e_L\}$, and $f\in Y_k^{\e,k'}$ then we can write
\begin{equation}\label{pp2}
\begin{split}
f^{\e,k'}&(\xi_1\e+\xi',\tau)=\eta^+_{[k'-1,k'+1]}(M)\\
&\times\frac{2^{-k'/2}\cdot\eta_{\leq k'-80}(\xi_1-M)}{\xi_1-M+i/2^{k'}}\int_\mathbb{R}e^{-iy_1\xi_1}h(y_1,\xi',\tau)\,dy_1+g,
\end{split}
\end{equation}
where $\xi_1,\tau\in\R$, $\xi'\in P_\e$, $h$ is supported in $\R\times S_{k,k'}^\e$,
\begin{equation}\label{pp22}
S_{k,k'}^\e=\{(\xi',\tau)\in P_\e\times\R:-\tau-|\xi'|^2\in[2^{2k'-10},2^{2k'+10}]\text{  and }|\xi'|\leq 2^{k+1}\},
\end{equation}
$M=M(\xi',\tau)=(-\tau-|\xi'|^2)^{1/2}$, and
\begin{equation}\label{pp3}
||g||_{X_k}+||h||_{L^1_{y_1}L^2_{\xi',\tau}}\leq (C/ \gamma_{k,k'})\cdot ||f^{\e,k'}||_{Y^{\e,k'}_k}.
\end{equation}
\end{Lemmar1}

\begin{proof}[Proof of Lemma \ref{Lemmar1}] Let
\begin{equation*}
h'(x,t)=2^{-k'/2}\mathcal{F}^{-1}[(\tau+|\xi|^2+i)\cdot f^{\e,k'}](x,t),
\end{equation*}
so
\begin{equation}\label{vc21}
\begin{cases}
&f^{\e,k'}(\xi_1\e+\xi',\tau)=\eta^+_{[k'-1,k'+1]}(\xi_1)\cdot \mathbf{1}_{[0,2^{k+1}]}(|\xi'|)\cdot \frac{2^{k'/2}}{\tau+|\xi|^2+i}\mathcal{F}(h')(\xi_1\e+\xi',\tau);\\
&\|h'\|_{L^{1,2}_\e}\leq (C/ \gamma_{k,k'})\cdot||f^{\e,k'}||_{Y^{\e,k'}_k}.
\end{cases}
\end{equation}
Let
\begin{equation*}
h''(y_1,\xi',\tau)=\int_{P_\e\times\R}h'(y_1\e+y',t)e^{-i(y'\cdot\xi'+t\tau)}\,dy'dt.
\end{equation*}
As in Lemma \ref{Lemmas1} (see \eqref{pp1}),
\begin{equation*}
||\eta_{\geq 2k'-79}(\tau+\xi^2)\cdot f^{\e,k'}||_{X_k}\leq (C/ \gamma_{k,k'})\cdot||f^{\e,k'}||_{Y^{\e,k'}_k}.
\end{equation*}
Thus it remains to write $\eta_{\leq2k'-80}(\tau+\xi^2)\cdot f^{\e,k'}$ as in \eqref{pp2}. Using \eqref{vc21}
\begin{equation}\label{vc22}
\begin{split}
&\eta_{\leq2k'-80}(\tau+\xi^2)\cdot f^{\e,k'}(\xi_1\e+\xi',\tau)\\
&=2^{k'/2}\cdot \eta^+_{[k'-1,k'+1]}(\xi_1)\cdot  \mathbf{1}_{[0,2^{k+1}]}(|\xi'|)\cdot \frac{\eta_{\leq2k'-80}(\tau+\xi^2)}{\tau+|\xi|^2+i}\int_{\R}h''(y_1,\xi',\tau)e^{-iy_1\xi_1}\,dy_1.
\end{split}
\end{equation}
Clearly, we may assume that $h''$ is supported in $\R\times S_{k,k'}^\e$. Let $M=(-\tau-|\xi'|^2)^{1/2}$ and approximate
\begin{equation}\label{vc23}
\begin{split}
&\eta^+_{[k'-1,k'+1]}(\xi_1)\cdot \mathbf{1}_{[0,2^{k+1}]}(|\xi'|)\cdot \frac{\eta_{\leq 2k'-80}(\xi_1^2-M^2)}{\xi_1^2-M^2+i}\\
&=\eta^+_{[k'-1,k'+1]}(M)\cdot \mathbf{1}_{[0,2^{k+1}]}(|\xi'|)\cdot \frac{\eta_{\leq k'-80}(\xi_1-M)}{\xi_1-M+i/2^{k'}}\cdot \frac{1}{2M}+E(\xi_1,\xi',\tau),
\end{split}
\end{equation}
where, with $\mu=|\xi_1^2-M^2|+1=\big|\tau+|\xi|^2\big|+1$,
\begin{equation*}
|E(\xi_1,\xi',\tau)|\leq  C\mathbf{1}_{[0,2^{k+1}]}(|\xi'|)\eta^+_{[k'-5,k'+5]}(\xi_1)\cdot \frac{\eta_{\leq 2k'}(\mu)}{\mu}\cdot \Big(\frac{\mu}{2^{2k'}}+\frac{1}{\mu}\Big).
\end{equation*}
We substitute \eqref{vc23} into \eqref{vc22} and notice that the error term corresponding to $E(\xi_1,\xi',\tau)$ can be bounded in $X_k$ (as  in Lemma \ref{Lemmas1}). The  main term in the right-hand side of \eqref{vc23} leads to the representation \eqref{pp2}, with $$h=\mathbf{1}_{[0,2^{k+1}]}(|\xi'|)\cdot 2^{k'}\cdot (2M)^{-1}\cdot h''.$$
\end{proof}

\section{Proof of Theorem \ref{Main2}}\label{proof}

For $\sigma\geq d/2$ we define the normed spaces
\begin{equation}\label{no5}
F^\sigma=\{u\in C(\mathbb{R}:\dot{B}^\infty):\|u\|_{F^\sigma}=\sum_{k\in\Z} \max(2^{dk/2},2^{\sigma k})\cdot \|\eta_k^{(d)}(\xi)\cdot \mathcal{F}(u)\|_{Z_k}<\infty\},
\end{equation}
and
\begin{equation}\label{no6}
\begin{split}
&N^\sigma=\{u\in C(\mathbb{R}:\dot{B}^\infty):\\
&\|u\|_{N^\sigma}=\sum_{k\in\Z} \max(2^{dk/2},2^{\sigma k})\cdot \|\eta_k^{(d)}(\xi)\cdot (\tau+|\xi|^2+i)^{-1}\cdot \mathcal{F}(u)\|_{Z_k}<\infty\}.
\end{split}
\end{equation}

For $\phi\in \dot{B}^\infty$ let $W(t)\phi\in
C(\mathbb{R}:\dot{B}^\infty)$ denote the solution of the free Schr\"{o}dinger evolution. Assume
$\psi:\mathbb{R}\to[0,1]$ is an even smooth function supported in
the interval $[-8/5,8/5]$ and equal to $1$ in the interval
$[-5/4,5/4]$. We prove first two linear estimates.

\newtheorem{Lemmaqq1}{Lemma}[section]
\begin{Lemmaqq1}\label{Lemmaqq1}
If $\sigma\geq d/2$ and $\phi\in \dot{B}^{\infty}$ then $\psi(t)\cdot [W(t)\phi]\in F^\sigma$ and
\begin{equation*}
\|\psi(t)\cdot [W(t)\phi]\|_{F^{\sigma}}\leq C\|\phi\|_{\dot{B}^\sigma}.
\end{equation*}
\end{Lemmaqq1}

\begin{proof}[Proof of Lemma \ref{Lemmaq1}] A straightforward computation shows that
\begin{equation*}
\mathcal{F}[\psi(t)\cdot (W(t)\phi)](\xi,\tau)=
\mathcal{F}_{(d)}(\phi)(\xi)\cdot \mathcal{F}_{(1)}(\psi)(\tau+|\xi|^2).
\end{equation*}
Then, directly from the definitions,
\begin{equation*}
\begin{split}
\|\psi(t)\cdot [W(t)\phi]\|_{F^{\sigma}}&=\sum_{k\in\mathbb{Z}}\max(2^{dk/2},2^{\sigma k})\|\eta_k^{(d)}(\xi)\mathcal{F}_{(d)}(\phi)(\xi)\mathcal{F}_{(1)}(\psi)(\tau+|\xi|^2)\|_{Z_k}\\
&\leq \sum_{k\in\mathbb{Z}}\max(2^{dk/2},2^{\sigma k})\|\eta_k^{(d)}(\xi)\mathcal{F}_{(d)}(\phi)(\xi)\mathcal{F}_{(1)}(\psi)(\tau+|\xi|^2)\|_{X_k}\\
&\leq C\sum_{k\in\mathbb{Z}}\max(2^{dk/2},2^{\sigma k})\|\eta_k^{(d)}(\xi)\cdot \mathcal{F}_{(d)}(\phi)(\xi)\|_{L^2}\\
&\leq C\|\phi\|_{\dot{B}^\sigma},
\end{split}
\end{equation*}
as desired.
\end{proof}

\newtheorem{Lemmaq3}[Lemmaqq1]{Lemma}
\begin{Lemmaq3}\label{Lemmaq3}
If $\sigma\geq 0$ and $u\in N^{\sigma}$ then $\psi(t)\cdot \int_0^tW(t-s)(u(s))\,ds\in F^\sigma$ and
\begin{equation*}
\Big|\Big|\psi(t)\cdot \int_0^tW(t-s)(u(s))\,ds\Big|\Big|_{F^{\sigma}}\leq C||u||_{N^{\sigma}}.
\end{equation*}
\end{Lemmaq3}

\begin{proof}[Proof of Lemma \ref{Lemmaq3}] A straightforward computation shows that
\begin{equation*}
\begin{split}
\mathcal{F}\Big[\psi(t)\cdot& \int_0^tW(t-s)(u(s))ds\Big](\xi,\tau)=\\
&c\int_\mathbb{R}\mathcal{F}(u)(\xi,\tau')\frac{\widehat{\psi}(\tau-\tau')-\widehat{\psi}(\tau+|\xi|^2)}{\tau'+|\xi|^2}d\tau',
\end{split}
\end{equation*}
where, for simplicity of notation, $\widehat{\psi}=\mathcal{F}_{(1)}(\psi)$. For $k\in\mathbb{Z}$ let
$$f_k(\xi,\tau')=\mathcal{F}(u)(\xi,\tau')\cdot \eta_k^{(d)}(\xi)\cdot (\tau'+|\xi|^2+i)^{-1}.$$
For $f\in Z_k$ let
\begin{equation}\label{ar202}
T(f)(\xi,\tau)=\int_\mathbb{R}f(\xi,\tau')\frac{\widehat{\psi}(\tau-\tau')-
\widehat{\psi}(\tau+|\xi|^2)}{\tau'+|\xi|^2}(\tau'+|\xi|^2+i)\,d\tau'.
\end{equation}
In view of the definitions, it suffices to prove that
\begin{equation}\label{ni5}
||T||_{Z_k\to Z_k}\leq C\text{ uniformly in }k\in\mathbb{Z}.
\end{equation}

To prove \eqref{ni5} we use the representation \eqref{mi1}. Assume first that $f=g_j$ is supported in $D_{k,j}$. Let $g_j^\#(\xi,\mu')=g_j(\xi,\mu'-|\xi|^2)$ and $[T(g)]^\#(\xi,\mu)=T(g)(\xi,\mu-|\xi|^2)$. Then,
\begin{equation}\label{ni6}
[T(g)]^\#(\xi,\mu)=\int_\mathbb{R}g_j^\#(\xi,\mu')\frac{\widehat{\psi}(\mu-\mu')-
\widehat{\psi}(\mu)}{\mu'}(\mu'+i)\,d\mu'.
\end{equation}
We use the elementary bound
\begin{equation*}
\Big|\frac{\widehat{\psi}(\mu-\mu')-\widehat{\psi}(\mu)}{\mu'}(\mu'+i)\Big|\leq C[(1+|\mu|)^{-4}+(1+|\mu-\mu'|)^{-4}].
\end{equation*}
Then, using \eqref{ni6},
\begin{equation*}
\begin{split}
|T(g)^\#(\xi,\mu)|&\leq C(1+|\mu|)^{-4}\cdot 2^{j/2}\Big[\int_{\mathbb{R}}|g_j^\#(\xi,\mu')|^2\,d\mu'\Big]^{1/2}\\
&+C\mathbf{1}_{[-2^{j+10},2^{j+10}]}(\mu)\int_{\mathbb{R}}|g_j^\#(\xi,\mu')|(1+|\mu-\mu'|)^{-4}\,d\mu'.
\end{split}
\end{equation*}
It follows from the definition of the spaces $X_k$ that
\begin{equation}\label{ni7}
||T||_{X_k\to X_k}\leq C\text{ uniformly in }k\in\mathbb{Z}_+,
\end{equation}
as desired.

Assume now that $f=f^{\mathbf{e}}\in Y^\mathbf{e}_k$,  $k\geq 100$, $\mathbf{e}\in\{\mathbf{e}_1,\ldots,\mathbf{e}_L\}$. We write
\begin{equation*}
f^{\mathbf{e}}(\xi,\tau')=\frac{\tau'+|\xi|^2}{\tau'+|\xi|^2+i}f^{\mathbf{e}}(\xi,\tau')+\frac{i}{\tau'+|\xi|^2+i}f^{\mathbf{e}}(\xi,\tau').
\end{equation*}
Using Lemma \ref{Lemmas1}, $||i(\tau'+|\xi|^2+i)^{-1}f^{\mathbf{e}}(\xi,\tau')||_{X_k}\leq C||f^\e||_{Y_k^\mathbf{e}}$. In view of \eqref{ar202} and \eqref{ni7}, for \eqref{ni5} it suffices to prove that
\begin{equation}\label{ni8}
\Big|\Big|\int_\mathbb{R}f^{\mathbf{e}}(\xi,\tau')\widehat{\psi}(\tau-\tau')\,d\tau'\Big|\Big|_{Z_k}+
\Big|\Big|\widehat{\psi}(\tau+|\xi|^2)\int_\mathbb{R}f^{\mathbf{e}}(\xi,\tau')\,d\tau'\Big|\Big|_{X_k}\leq C||f^{\mathbf{e}}||_{Y_k^{\mathbf{e}}}.
\end{equation}
The bound for the first term in the left-hand side of \eqref{ni8} follows easily from the definition. The bound for the second term in the left-hand side of \eqref{ni8} follows from \eqref{si1'} with $t=0$. 
\end{proof}

We prove now several nonlinear estimates. The main ingredients are the dyadic estimates in Lemma \ref{Lemmap1}, Lemma \ref{Lemmap2}, Lemma \ref{Lemmap3}, and Lemma \ref{Lemmaq1}. We reproduce these dyadic estimates below:

$\bullet$ if $k_1,k_2,k\in\Z$, $k_1\leq k_2+10$, $f_{k_1}\in Z_{k_1}$, and $f_{k_2}\in Z_{k_2}$, then
\begin{equation}\label{amy1}
2^{dk/2}\|\eta_k^{(d)}(\xi)\cdot (\widetilde{f}_{k_1}\ast f_{k_2})\|_{Z_k}\leq C2^{-|k_2-k|/4}(2^{dk_1/2}\|f_{k_1}\|_{Z_{k_1}})\cdot(2^{dk_2/2}\|f_{k_2}\|_{Z_{k_2}}),
\end{equation}
where $\mathcal{F}^{-1}(\widetilde{f}_{k_1})\in \{\mathcal{F}^{-1}(f_{k_1}),\overline{\mathcal{F}^{-1}(f_{k_1})}\}$.

$\bullet$ if $k_1,k_2,k\in\Z$, $k_1\leq k_2-10$, $|k-k_2|\leq 2$, $f_{k_1}\in Z_{k_1}$, and $f_{k_2}\in Z_{k_2}$ then
\begin{equation}\label{amy2}
\begin{split}
2^{dk/2}&\|\eta_k^{(d)}(\xi)\cdot (\tau+|\xi|^2+i)^{-1}\big[\widetilde{f}_{k_1}\ast [(\tau_2+|\xi_2|^2+i)f_{k_2}]\big]\|_{Z_k}\\
&\leq C(2^{dk_1/2}\|f_{k_1}\|_{Z_{k_1}})\cdot(2^{dk_2/2}\|f_{k_2}\|_{Z_{k_2}}),
\end{split}
\end{equation}
where $\mathcal{F}^{-1}(\widetilde{f}_{k_1})\in \{\mathcal{F}^{-1}(f_{k_1}),\overline{\mathcal{F}^{-1}(f_{k_1})}\}$.

$\bullet$ if $k_1,k_2,k\in\Z$, $k_1\leq k_2+10$, $f_{k_1}\in Z_{k_1}$, and $f_{k_2}\in Z_{k_2}$, then
\begin{equation}\label{amy3}
\begin{split}
2^{dk/2}&\|\eta_k^{(d)}(\xi)\cdot (\tau+|\xi|^2+i)^{-1}\big[ [(\tau_1+|\xi_1|^2+i)f_{k_1}]\ast f_{k_2}\big]\|_{Z_k}\\
&\leq C2^{-|k_2-k|/4}\cdot (2^{dk_1/2}\|f_{k_1}\|_{Z_{k_1}})\cdot(2^{dk_2/2}\|f_{k_2}\|_{Z_{k_2}}).
\end{split}
\end{equation}

$\bullet$ if $k_1,k_2,k_3,k\in\Z$, $f_{k_1}\in Z_{k_1},\,f_{k_2}\in Z_{k_2},\,f_{k_3}\in Z_{k_3}$, and
\begin{equation}\label{amy5}
\min(k,k_2,k_3)\leq k_1+20,
\end{equation}
then
\begin{equation}\label{amy4}
\begin{split}
&2^{k_2+k_3}\cdot 2^{dk/2}\|\eta_{k}^{(d)}(\xi)\cdot (\tau+|\xi|^2+i)^{-1}\cdot(\widetilde{f}_{k_1}\ast\widetilde{f}_{k_2}\ast \widetilde{f}_{k_3})\|_{Z_k}\\
&\leq C2^{-|\max(k_1,k_2,k_3)-k|/4}\cdot (2^{dk_1/2}\|f_{k_1}\|_{Z_{k_1}})\cdot(2^{dk_2/2}\|f_{k_2}\|_{Z_{k_2}})\cdot (2^{dk_3/2}\|f_{k_3}\|_{Z_{k_3}}),
\end{split}
\end{equation}
where $\mathcal{F}^{-1}(\widetilde{f}_{k_l})\in \{\mathcal{F}^{-1}(f_{k_l}),\overline{\mathcal{F}^{-1}(f_{k_l})}\}$, $l=1,2,3$.

For $\sigma\geq d/2$ let
\begin{equation*}
\overline{F}^{\sigma}=\{u\in C(\R:\dot{B}^{\infty}:\overline{u}\in F^\sigma \text{  and }||u||_{\overline{F}^\sigma}=||\overline{u}||_{F^\sigma}\}.
\end{equation*}

\newtheorem{Lemmaqq2}[Lemmaqq1]{Lemma}
\begin{Lemmaqq2}\label{Lemmaqq2}
(a) If $u,v\in F^{d/2}$ then $u\cdot v\in F^{d/2}$, and
\begin{equation*}
||u\cdot v||_{F^{d/2}}\leq C||u||_{F^{d/2}}\cdot ||v||_{F^{d/2}}.
\end{equation*}

(b) If $u,v\in F^{d/2}+\overline{F}^{d/2}$ then $u\cdot v\in F^{d/2}+\overline{F}^{d/2}$, and
\begin{equation*}
||u\cdot v||_{F^{d/2}+\overline{F}^{d/2}}\leq C||u||_{F^{d/2}+\overline{F}^{d/2}}\cdot ||v||_{F^{d/2}+\overline{F}^{d/2}}.
\end{equation*}

(c) If $u\in F^{d/2}+\overline{F}^{d/2}$ and $v,w\in F^{d/2}$ then $u\cdot2\sum_{l=1}^d\partial_{x_l}v\cdot \partial_{x_l}w\in N^{d/2}$ and
\begin{equation*}
\big|\big|u\cdot2\sum_{l=1}^d\partial_{x_l}v\cdot \partial_{x_l}w\big|\big|_{N^{d/2}}\leq C||u||_{F^{d/2}+\overline{F}^{d/2}}\cdot ||v||_{F^{d/2}}\cdot ||w||_{F^{d/2}}.
\end{equation*}
\end{Lemmaqq2}

\begin{proof}[Proof of Lemma \ref{Lemmaqq2}] For part (a), let $f_k=\eta_{k}^{(d)}(\xi)\cdot\mathcal{F}(u)$, $g_k=\eta_{k}^{(d)}(\xi)\cdot\mathcal{F}(v)$, $k\in\Z$. For part (a) it suffices to prove that for any $k_1,k_2\in\Z$
\begin{equation*}
\sum_{k\in\Z}2^{dk/2}||\eta_k^{(d)}(\xi)\cdot (f_{k_1}\ast g_{k_2})||_{Z_k}\leq C(2^{dk_1/2}||f_{k_1}||_{Z_{k_1}})\cdot (2^{dk_2/2}||g_{k_2}||_{Z_{k_2}}),
\end{equation*}
which follows easily from \eqref{amy1}. The proof  of part (b) is similar, using only \eqref{amy1} and the definitions.

For part (c), for $k\in\Z$ let $f_k=\eta_k(\xi)\cdot\mathcal{F}(u)$, $u_k=\mathcal{F}^{-1}(f_k)$, $g_k=\eta_k(\xi)\cdot\mathcal{F}(v)$, $v_k=\mathcal{F}^{-1}(g_k)$, $h_k=\eta_k(\xi)\cdot\mathcal{F}(w)$, $w_k=\mathcal{F}^{-1}(h_k)$. It suffices to prove that for any $k_1,k_2,k_3\in \Z$
\begin{equation}\label{jj1}
\begin{split}
\sum_{k\in\Z}2^{dk/2}&\big|\big|\eta_k^{(d)}(\xi)\cdot(\tau+|\xi|^2+i)^{-1}\cdot \mathcal{F}\big[\widetilde{u}_{k_1}\cdot 2\sum_{l=1}^d\partial_{x_l}v_{k_2}\cdot \partial_{x_l}w_{k_3}\big]\big|\big|_{Z_k}\\
&\leq C(2^{dk_1/2}||f_{k_1}||_{Z_{k_1}})\cdot (2^{dk_2/2}||g_{k_2}||_{Z_{k_2}})\cdot (2^{dk_3/2}||h_{k_3}||_{Z_{k_3}}),
\end{split}
\end{equation}
where $\widetilde{u}_{k_1}\in\{u_{k_1},\overline{u}_{k_1}\}$. If $\min(k_2,k_3)\leq k_1+20$ then \eqref{jj1} follows directly from \eqref{amy4}. Assume that
\begin{equation*}
\min(k_2,k_3)\geq k_1+20.
\end{equation*}
Using \eqref{amy4} again, we only need to bound the sum over $k\geq k_1+20$. In this case we use the identity
\begin{equation}\label{null}
2\sum_{l=1}^d\partial_{x_l}v_{k_2}\cdot \partial_{x_l}w_{k_3}=H(v_{k_2}\cdot w_{k_3})-w_{k_3}\cdot Hv_{k_2}-v_{k_2}\cdot Hw_{k_3},
\end{equation}
where $H=i\partial_t+\Delta_x$. We estimate the sum over $k\geq k_1+20$ corresponding to the term $H(v_{k_2}\cdot w_{k_3})$ using \eqref{amy1} and \eqref{amy2}. We estimate the sums over $k\geq k_1+20$ corresponding to the terms $w_{k_3}\cdot Hv_{k_2}$ and $v_{k_2}\cdot Hw_{k_3}$ using \eqref{amy1}, \eqref{amy2}, and \eqref{amy3}. The bound \eqref{jj1} follows easily.
\end{proof}

Let
\begin{equation*}
\mathcal{N}(u)=2\overline{u}(1+u\overline{u})^{-1}\sum_{j=1}^d(\partial_{x_j}u)^2
\end{equation*}
denote the nonlinear term in \eqref{Sch5}. It follows from Lemma \ref{Lemmaqq2} that
\begin{equation}\label{jj2}
||\mathcal{N}(u)-\mathcal{N}(v)||_{N^{d/2}}\leq C\epsilon^2||u-v||_{F^{d/2}}
\end{equation}
for any $u,v\in B_{F^{d/2}}(0,\epsilon)$, $\epsilon\ll 1$, and
\begin{equation}\label{jj3}
||\partial_{x_l}^m\mathcal{N}(u)||_{N^{d/2}}\leq C\epsilon^2||\partial_{x_l}^mu||_{F^{d/2}}+C(m,||u||_{F^{d/2+m-1}}),
\end{equation}
for any $l\in\{1,\ldots,d\}$, $m\in 1,2,\ldots$, and  $u\in B_{F^{d/2}}(0,\epsilon)\cap F^{d/2+m}$. The bounds  \eqref{jj2} and \eqref{jj3}, together with the imbedding $F^{d/2}\hookrightarrow C(\R:\dot{B}^{d/2})$ (which follows from Lemma \ref{Lemmaa2}) are sufficient to construct the solution $u\in C(\R:\dot{B}^\infty)$ in Theorem \ref{Main2} and prove the Lipschitz bound \eqref{by1} (see, for example, \cite[Section 5]{IoKe2} for the standard recursive argument).

The uniqueness of solutions in $C(\R:\dot{B}^\infty\cap B_{\dot{B}^{d/2}}(0,\overline{\epsilon}_1))$, for $\overline{\epsilon}_1$ sufficiently small, follows from the following simple observation: if $u\in C([-T,T]:\dot{B}^\infty\cap B_{\dot{B}^{d/2}}(0,\overline{\epsilon}_1))$ is a solution of the equation \eqref{Sch5} then there is $\delta=\delta(||u||_{L^\infty_t\dot{B}^{d/2+100}})$ with the  property that
\begin{equation}\label{jj4}
||u||_{F^{d/2}[t_0-\delta,t_0+\delta]}\leq C\overline{\epsilon}_1\text{ for any }t_0\in [-T+\delta, T-\delta].
\end{equation}
See, for example, \cite[Section 10]{IoKe} for such an argument. The  uniqueness of solutions then follows from \eqref{jj4} and \eqref{jj2}.

\section{Maximal function and local smoothing estimates}\label{maxsmo}

In this section we  prove two lemmas that will be used extensively in the bi\-li\-near and the trilinear estimates in the following three sections. For $l=1,\ldots,d$ and $k\in \Z_+$ let $\Xi_k^{(l)}=2^k\cdot\Z^l$. Let $\chi^{(1)}:\mathbb{R}\to[0,1]$ denote a fixed smooth function supported in the interval $[-2/3,2/3]$ with the property that
\begin{equation*}
\sum_{n\in\Z}\chi^{(1)}(\xi-n)\equiv 1\text{ on }\R.
\end{equation*}
Let $\widetilde{\chi}^{(1)}:\mathbb{R}\to[0,1]$ denote a fixed smooth function supported in the interval $[-4,4]$ and equal to $1$ in the interval $[-3,3]$. Let $\chi^{(l)},\widetilde{\chi}^{(l)}:\mathbb{R}^l\to[0,1]$,
\begin{equation}\label{nn7}
\chi^{(l)}(\xi)=\chi^{(1)}(\xi_1)\cdot\ldots\cdot \chi^{(1)}(\xi_l)\text{ and }\widetilde{\chi}^{(l)}(\xi)=\widetilde{\chi}^{(1)}(\xi_1)\cdot\ldots\cdot \widetilde{\chi}^{(1)}(\xi_l).
\end{equation}
For $k\in \Z_+$ and $n\in\Xi_k^{(l)}$ we define
\begin{equation*}
\chi_{k,n}^{(l)}(\xi)=\chi^{(l)}((\xi-n)/2^k)\text{ and }\widetilde{\chi}^{(l)}_{k,n}(\xi)=\widetilde{\chi}^{(l)}((\xi-n)/2^k).
\end{equation*}
Clearly,
\begin{equation*}
\sum_{n\in\Xi_k}\chi^{(l)}_{k,n}\equiv 1\text{ on }\R^l.
\end{equation*}
For simplicity of notation, we let $\chi_{k,n}=\chi_{k,n}^{(d)}$
and $\Xi_k=\Xi_k^{(d)}$.

We start with a maximal function estimate.

\newtheorem{Lemmaa1}{Lemma}[section]
\begin{Lemmaa1}\label{Lemmaa1}
If $d\geq 3$, $k\in\Z$, $f\in Z_k$, and $\mathbf{e}'\in\mathbb{S}^{d-1}$ then
\begin{equation}\label{pr40}
||\mathcal{F}^{-1}(f)||_{L^{2,\infty}_{\mathbf{e}'}}\leq
C2^{(d-1)k/2}\|f\|_{Z_k}.
\end{equation}
In addition, if $k\geq 100$, $k_1\in[0,k+10d]\cap\Z$, and $f\in X_k$ then
\begin{equation}\label{mm1}
\big[\sum_{n\in\Xi_{k_1}}||\mathcal{F}^{-1}(\chi_{k_1,n}(\xi)\cdot f)||_{L^{2,\infty}_{\e'}}^2\big]^{1/2}\leq C2^{(d-1)k_1/2}\cdot 2^{(k-k_1)/2}\|f\|_{X_k}.
\end{equation}
If $k\geq 100$, $k_1\in[0,k+10d]\cap\Z$, and $f\in Z_k$ then
\begin{equation}\label{mm1.1}
\big[\sum_{n\in\Xi_{k_1}}||\mathcal{F}^{-1}[\chi_{k_1,n}(\xi)\cdot f\cdot\eta_{\leq k+k_1}(\tau+|\xi|^2)]||_{L^{2,\infty}_{\e'}}^2\big]^{1/2}\leq C2^{(d-1)k_1/2}\cdot 2^{(k-k_1)/2}\|f\|_{Z_k}.
\end{equation}
\end{Lemmaa1}

\begin{proof}[Proof of Lemma \ref{Lemmaa1}] We use the representation \eqref{mi1}
and assume first that $f=g_j$ is supported in $D_{k,j}$. For \eqref{pr40},
it suffices to prove that
\begin{equation}\label{pr41}
|| \mathcal{F}^{-1}(g_j)||_{L^{2,\infty}_{\mathbf{e}'}}\leq
C2^{(d-1)k/2}\cdot 2^{j/2}\|g_j\|_{L^2}.
\end{equation}
We define $g_j^\#(\xi,\mu)=g_j(\xi,\mu-|\xi|^2)$. The left-hand side of \eqref{pr41} is dominated by
\begin{equation*}
\begin{split}
\int_{[-2^{j+1},2^{j+1}]}\Big|\Big|\int_{\mathbb{R}^d}g^\#_{j}(\xi,\mu)e^{ix\cdot \xi}e^{-it|\xi|^2}\,d\xi\Big|\Big|_{L^{2,\infty}_{\mathbf{e}'}}\,d\mu.
\end{split}
\end{equation*}
Thus, for \eqref{pr41} it suffices to prove that
\begin{equation}\label{ar400}
\Big|\Big|\int_{\mathbb{R}^d}h(\xi)e^{ix\cdot\xi}e^{-it|\xi|^2}\,d\xi\Big|\Big|_{L^{2,\infty}_{\mathbf{e}'}}\leq C2^{(d-1)k/2}\cdot ||h||_{L^2_\xi},
\end{equation}
for any function $h$ supported in the set $\{\xi\in\mathbb{R}^d:|\xi|\leq 2^{k+1}\}$.

To prove \eqref{ar400}, using a standard $TT^\ast$ argument, it suffices to show that
\begin{equation}\label{ar401}
\begin{split}
\Big|&\Big|\int_{\mathbb{R}^{d-1}\times\mathbb{R}}e^{ix_1\xi_1}e^{ix'\cdot \xi'}e^{-it(\xi_1^2+|\xi'|^2)}\\
&\times \eta_0(\xi_1/2^{k+1})\cdot \eta_0( |\xi'|/2^{k+1})\,d\xi_1d\xi'\Big|\Big|_{L^1_{x_1}L^\infty_{x',t}}\leq C2^{(d-1)k}.
\end{split}
\end{equation}
By stationary phase (one  may also rescale to $k=0$), for any $x'\in\mathbb{R}^{d-1}$ and $x_1\in\R$
\begin{equation*}
\Big|\int_{\mathbb{R}^{d-1}}e^{ix'\cdot \xi'}e^{-it|\xi'|^2}\eta_0( |\xi'|/2^{k+1})\,d\xi'\Big|\leq C\min(2^{(d-1)k},|t|^{-(d-1)/2}),
\end{equation*}
and
\begin{equation*}
\Big|\int_{\mathbb{R}}e^{ix_1\cdot \xi_1}e^{-it\xi_1^2}\eta_0(\xi_1/2^{k+1})\,d\xi_1\Big|\leq C\min(2^{k},|t|^{-1/2}).
\end{equation*}
In addition, by integration by parts, if $|x_1|\geq 2^{k+10}|t|$ then
\begin{equation*}
\Big|\int_{\mathbb{R}}e^{ix_1\cdot \xi_1}e^{-it\xi_1^2}\eta_0(\xi_1/2^{k+1})\,d\xi_1\Big|\leq C2^k(1+2^k|x_1| )^{-2}.
\end{equation*}
Let $K(x_1,x',t)$ denote the function in the left-hand side of \eqref{ar401}. In view of the three bounds above,
\begin{equation*}
\sup_{t\in\R,\,x'\in\mathbb{R}^{d-1}}|K(x_1,x',t)|\leq C2^{dk}(1+2^k|x_1| )^{-2}+C2^{dk/2}|x_1|^{-d/2}\cdot\mathbf{1}_{[2^{-k},\infty)}( |x_1| ).
\end{equation*}
The bound \eqref{ar401} follows since $d\geq 3$.

We prove now the bound \eqref{mm1}, assuming $k\geq 100$ and $f=g_j$ is supported in $D_{k,j}$. Clearly we may assume $k_1\leq k-10d$, and it suffices to prove that for any $n\in\Xi_{k_1}$,
\begin{equation}\label{mm2}
|| \mathcal{F}^{-1}(\chi_{k_1,n}(\xi)\cdot g_j)||_{L^{2,\infty}_{\mathbf{e}'}}\leq
C2^{(d-1)k_1/2}2^{(k-k_1)/2}\cdot 2^{j/2}\|g_j\|_{L^2}.
\end{equation}
By the same argument as before, for \eqref{mm2} it suffices to prove that
\begin{equation}\label{mm3}
\Big|\Big|\int_{\mathbb{R}^d}h(\xi)e^{ix\cdot\xi}e^{-it|\xi+n|^2}\,d\xi\Big|\Big|_{L^{2,\infty}_{\mathbf{e}'}}\leq C2^{(d-1)k_1/2}2^{(k-k_1)/2}\cdot ||h||_{L^2_\xi},
\end{equation}
for any function $h$ supported in the set $\{\xi\in\mathbb{R}^d:|\xi|\leq 2^{k_1}\}$ and any vector $n\in \R^d$ with $|n|\leq 2^{k+2}$. As before, for \eqref{mm3}, it suffices to show that
\begin{equation}\label{mm4}
\begin{split}
\Big|&\Big|\int_{\mathbb{R}^{d-1}\times\mathbb{R}}e^{ix_1\xi_1}e^{ix'\cdot \xi'}e^{-it(\xi_1^2+|\xi'|^2)}e^{-2itn_1\xi_1}e^{-2itn'\cdot\xi'}\\
&\times \eta_0(\xi_1/2^{k_1})\cdot \eta_0( |\xi'|/2^{k_1})\,d\xi_1d\xi'\Big|\Big|_{L^1_{x_1}L^\infty_{x',t}}\leq C2^{(d-1)k_1}2^{k-k_1}.
\end{split}
\end{equation}
By stationary phase, for any $x'\in\mathbb{R}^{d-1}$ and $x_1\in\R$,
\begin{equation*}
\Big|\int_{\mathbb{R}^{d-1}}e^{ix'\cdot \xi'}e^{-2itn'\cdot\xi'}e^{-it|\xi'|^2}\eta_0( |\xi'|/2^{k_1})\,d\xi'\Big|\leq C\min(2^{(d-1)k_1},|t|^{-(d-1)/2}),
\end{equation*}
and
\begin{equation*}
\Big|\int_{\mathbb{R}}e^{ix_1\xi_1}e^{-2itn_1\xi_1}e^{-it\xi_1^2}\eta_0(\xi_1/2^{k_1})\,d\xi_1\Big|\leq C\min(2^{k_1},|t|^{-1/2}).
\end{equation*}
In addition, by integration by parts, if $|x_1|\geq 2^{k+10}|t|$ then
\begin{equation*}
\Big|\int_{\mathbb{R}}e^{ix_1\cdot \xi_1}e^{-2itn_1\xi_1}e^{-it\xi_1^2}\eta_0(\xi_1/2^{k_1})\,d\xi_1\Big|\leq C2^{k_1}(1+2^{k_1}|x_1| )^{-d/2}.
\end{equation*}
Let $K_1(x_1,x',t)$ denote the function in the left-hand side of \eqref{mm4}. In view of the three bounds above,
\begin{equation*}
\sup_{t\in\R,\,x'\in\mathbb{R}^{d-1}}|K_1(x_1,x',t)|\leq C2^{dk_1}\cdot\mathbf{1}_{[0,2^{k-2k_1}]}(|x_1|)+C2^{dk/2}|x_1|^{-d/2}\cdot\mathbf{1}_{[2^{k-2k_1},\infty)}( |x_1| ).
\end{equation*}
The bound \eqref{mm4} follows since $d\geq 3$.

Assume now that $f=f^{\e,k'}\in Y_k^{\e,k'}$, $k\geq 100$, $k'\in T_k$, $\e\in\{\e_1,\ldots,\e_L\}$. It suffices to prove the stronger bound \eqref{mm1.1}, and we may assume $k_1\leq [30d,k-30d]\cap\Z$. We assume first
\begin{equation}\label{ff1}
k_1\leq k'.
\end{equation}
We fix an arbitrary orthonormal basis in $P_\e$ and use it to define an isomorphism
\begin{equation*}
\Phi^\e:\R^{d-1}\to P_\e.
\end{equation*}
For $k\in \Z_+$ let $\Xi_k^\e=\Phi^\e(\Xi_k^{(d-1)})\subseteq P_\e$. For $n'\in\Xi_k^\e$ let
\begin{equation}\label{nn6}
\chi_{k,n'}^\e=\chi_{k,(\Phi^\e)^{-1}(n')}^{(d-1)}\circ(\Phi^\e)^{-1}\text{ and }\widetilde{\chi}_{k,n}^\e=\widetilde{\chi}_{k,(\Phi^\e)^{-1}(n')}^{(d-1)}\circ(\Phi^\e)^{-1},
\end{equation}
$\chi_{k,n'}^\e,\widetilde{\chi}_{k,n'}^\e:P_\e\to[0,1]$ (compare with the notation at the beginning of the section). We write $\xi=\xi_1\e+\xi'$, $x=x_1\e+x'$, $\xi_1,x_1\in\R$. $\xi',x'\in P_\e$. For \eqref{mm1.1} it suffices to prove that
\begin{equation}\label{vy1}
\begin{split}
\Big[&\sum_{n_1\in\Xi^{(1)}_{k_1},\,n_1\in[2^{k'-10},2^{k'+10}]}\,\,\,\sum_{n'\in\Xi^\e_{k_1},\,|n'|\leq 2^{k+10}}||\mathcal{F}^{-1}[\chi^\e_{k_1,n'}(\xi')\cdot\chi^{(1)}_{k_1,n_1}(\xi_1)\\
&\times f^{\e,k'}\cdot\eta_{\leq k+k_1}(\tau+|\xi|^2)]||_{L^{2,\infty}_{\e'}}^2\Big]^{1/2}\leq C2^{(d-1)k_1/2}2^{(k-k_1)/2}||f^{\e,k'}||_{Y_k^{\e,k'}}.
\end{split}
\end{equation}
We use Lemma \ref{Lemmas1} and Lemma \ref{Lemmar1} and notice that we  may replace $\eta_{\leq k+k_1}(\xi_1^2-M^2)\cdot \eta_{\leq k'-80}(\xi_1-M)$ by $\eta_{\leq k_1-50}(\xi_1-M)$, at the expense of $C(k-k'+10)$ error terms in $X_k$. The contributions of these error terms are controlled using \eqref{mm1} and the large factor $\gamma_{k,k'}$ in \eqref{v2}. For simplicity  of notation, in the rest of this  proof we let $\sum_{n_1,n'}$ denote the sum over $n_1$ and $n'$ as in \eqref{vy1}.  Using Lemma \ref{Lemmar1}, for \eqref{vy1} it suffices to prove the stronger bound
\begin{equation}\label{vy2}
\begin{split}
\Big[\sum_{n_1,n'}||\mathcal{F}^{-1}[\chi^\e_{k_1,n'}(\xi')\cdot\chi^{(1)}_{k_1,n_1}(\xi_1)\cdot \frac{2^{-k'/2}\cdot\eta_{\leq k_1-50}(\xi_1-M)}{\xi_1-M+1/2^{k'}}\cdot h'(\xi',\tau)]||_{L^{2,\infty}_{\e'}}^2\Big]^{1/2}\\
\leq C2^{(d-1)k_1/2}2^{(k-k_1)/2}||h'||_{L^2},
\end{split}
\end{equation}
for any $h'\in L^2(P_\e\times\R)$. We notice that $\chi^{(1)}_{k_1,n_1}(\xi_1)\cdot \eta_{\leq k_1-50}(\xi_1-M)=\chi^{(1)}_{k_1,n_1}(\xi_1)\cdot \eta_{\leq k_1-50}(\xi_1-M)\cdot\widetilde{\chi}^{(1)}_{k_1,n_1}(M)$. Thus the left-hand side of \eqref{vy2} is equal to
\begin{equation*}
\begin{split}
C\Big[\sum_{n_1,n'}\Big|\Big|\int_{R\times P_\e\times\R}e^{ix_1\xi_1}&e^{ix'\cdot\xi'}e^{it\tau}\chi^\e_{k_1,n'}(\xi')\cdot\chi^{(1)}_{k_1,n_1}(\xi_1)\cdot \widetilde{\chi}^{(1)}_{k_1,n_1}(M)\\
&\times\frac{2^{-k'/2}\cdot\eta_{\leq k_1-50}(\xi_1-M)}{\xi_1-M+1/2^{k'}}\cdot h'(\xi',\tau)\,d\xi_1d\xi' d\tau\Big|\Big|_{L^{2,\infty}_{\e'}}^2\Big]^{1/2}.
\end{split}
\end{equation*}
We may disregard the factor $\chi^{(1)}_{k_1,n_1}(\xi_1)$ in the integral above, and  integrate the variable $\xi_1$ first. The  left-hand side of  \eqref{vy2} is dominated by
\begin{equation*}
\begin{split}
C\Big[\sum_{n_1,n'}\Big|\Big|\int_{P_\e\times\R}e^{ix_1M}&e^{ix'\cdot\xi'}e^{it\tau}\chi^\e_{k_1,n'}(\xi')\cdot \widetilde{\chi}^{(1)}_{k_1,n_1}(M)\cdot 2^{-k'/2}h'(\xi',\tau)\,d\xi' d\tau\Big|\Big|_{L^{2,\infty}_{\e'}}^2\Big]^{1/2}.
\end{split}
\end{equation*}
We make the substitutions $\tau=-\mu^2-|\xi'|^2$ (so $M=\mu$) and $h''(\xi,\mu)=2^{-k'/2}\cdot\mu\cdot h'(\xi',-\mu^2-|\xi'|^2)$ (so $||h''||_{L^2}\approx ||h'||_{L^2}$). For \eqref{vy2} it suffices to prove that
\begin{equation*}
\begin{split}
\Big[\sum_{n_1,n'}\Big|\Big|\int_{P_\e\times\R}e^{ix_1\mu}e^{ix'\cdot\xi'}e^{-it(\mu^2+|\xi'|^2)}\chi^\e_{k_1,n'}(\xi')\cdot \widetilde{\chi}^{(1)}_{k_1,n_1}(\mu)\cdot h''(\xi',\mu)\,d\xi' d\mu\Big|\Big|_{L^{2,\infty}_{\e'}}^2\Big]^{1/2}\\
\leq C2^{(d-1)k_1/2}2^{(k-k_1)/2}||h''||_{L^2},
\end{split}
\end{equation*}
which follows from \eqref{mm3}.

If $k'\leq k_1$ then, using the large factor $\gamma_{k,k'}$ in \eqref{v2}, it suffices to prove the bound \eqref{mm1.1} in the case $k_1=k'$. This was already proved before.
\end{proof}

We prove now a local smoothing estimate.

\newtheorem{Lemmas2}[Lemmaa1]{Lemma}
\begin{Lemmas2}\label{Lemmas2}
If $k\in \Z$, $\mathbf{e}'\in\mathbb{S}^{d-1}$,
$l\in[-1,40]\cap\Z$, and $f\in Z_k$ then
\begin{equation}\label{gu40main}
\|\mathcal{F}^{-1}[f\cdot
\eta_{1}(\xi\cdot\mathbf{e}'/2^{k-l})]\|_{L^{\infty,2}_{\mathbf{e}'}}\leq
C2^{-k/2}\|f\|_{Z_k}.
\end{equation}
\end{Lemmas2}
\begin{proof}[Proof of Lemma \ref{Lemmas2}] Since $L^{p,q}_{\e'}\equiv L^{p,q}_{-\e'}$, for \eqref{gu40main} it suffices to prove that
\begin{equation}\label{gu40}
\|\mathcal{F}^{-1}[f\cdot
\eta_{1}^+(\xi\cdot\mathbf{e}'/2^{k-l})]\|_{L^{\infty,2}_{\mathbf{e}'}}\leq
C2^{-k/2}\|f\|_{Z_k}.
\end{equation}
We use the representation \eqref{mi1}. Assume first that $f=g_j$.
In view of the definitions, it suffices to prove that if $j\geq 0$
and $g_j$ is supported in $D_{k,j}$ then
\begin{equation}\label{gu41}
\Big|\Big|\int_{\mathbb{R}^d}e^{ix\cdot\xi}e^{it\tau}
g_j(\xi,\tau)\cdot\eta_{1}^+(\xi\cdot\mathbf{e}'/2^{k-l}) \,d\xi
d\tau\Big|\Big|_{L^{\infty,2}_{\e'}}\leq
C2^{-k/2}2^{j/2}\|g_j\|_{L^2}.
\end{equation}
Let $g_j^\#(\xi,\mu)=g_j(\xi,\mu-|\xi|^2)$. By H\"{o}lder's
ine\-qua\-lity, for \eqref{gu41} it suffices to prove that
\begin{equation}\label{vy10}
\Big|\Big|\int_{\mathbb{R}^d}e^{ix\cdot\xi}e^{-it|\xi|^2}h(\xi)
\cdot\eta^+_1(\xi\cdot\e'/2^{k-l})\,d\xi
d\tau\Big|\Big|_{L^{\infty,2}_{\e'}}\leq C2^{-k/2}\|h\|_{L^2},
\end{equation}
which follows easily using Plancherel theorem and a change of variables.

Assume now that $k\leq 100$, $f=f^{\mathbf{e},k'}\in
Y^{\mathbf{e},k'}_k$, $k'\in T_k$,
$\mathbf{e}\in\{\mathbf{e}_1,\ldots,\mathbf{e}_L\}$. The estimates
in Lemma \ref{Lemmar1} show that
\begin{equation*}
||f^{\mathbf{e},k'}\cdot [\eta^+_{[k-50,k+1]}(\xi\cdot\mathbf{e}')-\eta^+_{[k-50,k+1]}((M\e+\xi')\cdot\mathbf{e}')]||_{X_k}\leq C||f^{\mathbf{e},k'}||_{Y_k^{\mathbf{e},k'}}.
\end{equation*}
Since \eqref{gu40} was already proved for $f\in X_k$, it suffices to show that
\begin{equation}\label{gu47}
\begin{split}
\Big|\Big|&\int_{\R\times P_\e\times\mathbb{R}}e^{ix_1\xi_1}e^{ix'\cdot\xi'}e^{it\tau}f^{\mathbf{e},k'}(\xi_1\mathbf{e}+\xi',\tau)\\
&\times\eta^+_{[k-50,k+1]}((M\e+\xi')\cdot\e')\,d\xi_1d\xi' d\tau\Big|\Big|_{L^{\infty,2}_{\e'}}\leq C2^{-k/2}\|f^{\mathbf{e},k'}\|_{Y_k^{\mathbf{e},k'}}.
\end{split}
\end{equation}
We use now the representation in Lemma \ref{Lemmar1}, and  integrate the variable $\xi_1$ first in the left-hand side of \eqref{gu47}. For \eqref{gu47} it suffices to prove the stronger bound
\begin{equation}\label{gu48}
\begin{split}
\Big|\Big|\int_{P_\e\times\mathbb{R}}&e^{ix_1M}e^{ix'\cdot\xi'}e^{it\tau}\cdot \eta^+_{[k'-1,k'+1]}(M)\cdot \eta_0(|\xi'|/2^{k+2})\cdot 2^{-k'/2}h'(\xi',\tau)\\
&\times\eta^+_{[k-50,k+1]}((M\e+\xi')\cdot\e')\,d\xi' d\tau\Big|\Big|_{L^{\infty,2}_{\e'}}\leq C2^{-k/2}\|h'\|_{L^2_{\xi,\tau}},
\end{split}
\end{equation}
for any $h'\in L^2(P_\e\times\R)$. We  may the substitutions $\tau=-\mu^2-|\xi'|^2$ (so $M=\mu$), and $h''(\xi',\mu)=2^{-k'/2}\cdot\mu\cdot h'(\xi',-\mu^2-|\xi'|^2)$ (so $||h''||_{L^2}\approx ||h'||_{L^2}$). For \eqref{gu48} it suffices to prove that
\begin{equation*}
\begin{split}
\Big|\Big|\int_{P_\e\times\mathbb{R}}&e^{ix_1\mu}e^{ix'\cdot\xi'}e^{-it(\mu^2+|\xi'|^2)}\cdot \eta^+_{[k'-1,k'+1]}(\mu)\cdot \eta_0(|\xi'|/2^{k+2})\cdot h''(\xi',\mu)\\
&\times\eta^+_{[k-50,k+1]}((\mu\e+\xi')\cdot\e')\,d\xi' d\mu\Big|\Big|_{L^{\infty,2}_{\e'}}\leq C2^{-k/2}\|h''\|_{L^2_{\xi,\tau}},
\end{split}
\end{equation*}
which follows from \eqref{vy10}.
\end{proof}

\section{Dyadic bilinear estimates, I}\label{bilin1}

In this section we prove several dyadic bilinear estimates. We assume in the rest of this section that $d\geq 3$. We record first a simple $L^2$ estimate (see, for example, Lemma 6.1 (a) in \cite{IoKe} for the proof): if $k_1,k_2,k\in\Z$, $j_1,j_2,j\in\Z_+$, and $g_{k_1,j_1}$, $g_{k_2,,j_2}$ are $L^2$ functions supported in $D_{k_1,j_1}$ and $D_{k_2,j_2}$,  then
\begin{equation}\label{bt0}
||\mathbf{1}_{D_{k,j}}\cdot (\widetilde{g}_{k_1,j_1}\ast \widetilde{g}_{k_2,j_2})||_{L^2}\leq C2^{d\cdot \min(k_1,k_2,k)/2}2^{\min(j_1,j_2,j)/2}||g_{k_1,j_1}||_{L^2}\cdot||g_{k_2,j_2}||_{L^2},
\end{equation}
where $\mathcal{F}^{-1}(\widetilde{g}_{k_l,j_l})\in\{\mathcal{F}^{-1}(g_{k_l,j_l}),\overline{\mathcal{F}^{-1}(g_{k_l,j_l})}\}$, $l=1,2$.

For any $k\in\Z$, $j\in\Z_+$, and $f_k\in Z_k$ we let
\begin{equation}\label{qq9}
f_{k,\leq j}(\xi,\tau)=f_k(\xi,\tau)\cdot \eta_{\leq j}(\tau+|\xi|^2)\text{ and }f_{k,\geq j}(\xi,\tau)=f_k(\xi,\tau)\cdot \eta_{\geq j}(\tau+|\xi|^2).
\end{equation}
We will often use the following simple estimate.

\newtheorem{Lemmap0}{Lemma}[section]
\begin{Lemmap0}\label{Lemmap0}
If $k_1,k_2\in \Z$, $k_1\leq k_2+C$, $j_1,j_2\in \Z_+$, $f_{k_1}\in Z_{k_1}$, and $f_{k_2}\in Z_{k_2}$ then
\begin{equation}\label{qq1}
\begin{split}
&\|\widetilde{f}_{k_1,\geq j_1}\ast \widetilde{f}_{k_2,\geq j_2})\|_{L^2}\\
&\leq C(2^{j_2/2}+2^{(k_1+k_2)/2})^{-1}(\beta_{k_1,j_1}\cdot\beta_{k_2,j_2})^{-1}\cdot (2^{dk_1/2}\|f_{k_1}\|_{Z_{k_1}})\cdot(\|f_{k_2}\|_{Z_{k_2}}),
\end{split}
\end{equation}
where $\mathcal{F}^{-1}(\widetilde{f}_{k_l,\geq j_l})\in \{\mathcal{F}^{-1}(f_{k_l,\geq j_1}),\overline{\mathcal{F}^{-1}(f_{k_l,\geq j_l})}\}$, $l=1,2$.
\end{Lemmap0}

\begin{proof}[Proof of Lemma \ref{Lemmap0}] If $k_2\geq 100$ then, in view of \eqref{mi2}, we may assume that
\begin{equation}\label{qq2}
f_{k_2}\text{ is supported in }I_{k_2}^{(d)}\times\R\cap \{(\xi_2,\tau_2):|\xi_2-v|\leq 2^{k_2-50}\}\text{ for some }v\in I_{k_2}^{(d)}.
\end{equation}
Let $\widehat{v}=v/|v|$. Then, for $k_2\geq 100$, using Lemma \ref{Lemmas2} (and \eqref{gu41} when $j_2\geq 2k_2$),
\begin{equation}\label{qq3}
\|\mathcal{F}^{-1}(\widetilde{f}_{k_2,\geq j_2})\|_{L^{\infty,2}_{\widehat{v}}}\leq C2^{-k_2/2}\beta_{k_2,j_2}^{-1}\|f_{k_2}\|_{Z_{k_2}}.
\end{equation}
Using Lemma \ref{Lemmaa1} (and \eqref{pr41} when $j_1\geq 2k_1$),
\begin{equation}\label{qq4}
\|\mathcal{F}^{-1}(\widetilde{f}_{k_1,\geq j_1})\|_{L^{2,\infty}_{\widehat{v}}}\leq C2^{(d-1)k_1/2}\beta_{k_1,j_1}^{-1}\|f_{k_1}\|_{Z_{k_1}}.
\end{equation}
Using the definition,
\begin{equation}\label{qq5}
\|\mathcal{F}^{-1}(\widetilde{f}_{k_2,\geq j_2})\|_{L^{2}}\leq C2^{-j_2/2}\beta_{k_2,j_2}^{-1}\|f_{k_2}\|_{Z_{k_2}}.
\end{equation}
Finally, using Lemma \ref{Lemmaa2} (and \eqref{si1.2} when $j_1\geq 2k_1$),
\begin{equation}\label{qq6}
\|\mathcal{F}^{-1}(\widetilde{f}_{k_1,\geq j_1})\|_{L^{\infty}}\leq C2^{dk_1/2}\beta_{k_1,j_1}^{-1}\|f_{k_1}\|_{Z_{k_1}}.
\end{equation}
The bound \eqref{qq1} follows by using \eqref{qq3} and \eqref{qq4} when $k_1+k_2\geq j_2$, and \eqref{qq5} and \eqref{qq6} when $k_1+k_2\leq j_2$ (if $k_2\leq 100$ we always use \eqref{qq5} and \eqref{qq6}).
\end{proof}

Our next bilinear estimate is the main ingredient in the proof of the algebra properties in Lemma \ref{Lemmaqq2}.

\newtheorem{Lemmap1}[Lemmap0]{Lemma}
\begin{Lemmap1}\label{Lemmap1}
If $k_1,k_2,k\in\Z$, $k_1\leq k_2+10$, $f_{k_1}\in Z_{k_1}$, and $f_{k_2}\in Z_{k_2}$ then
\begin{equation}\label{bt1main}
2^{dk/2}\|\eta_k^{(d)}(\xi)\cdot (\widetilde{f}_{k_1}\ast f_{k_2})\|_{Z_k}\leq C2^{-|k_2-k|/4}(2^{dk_1/2}\|f_{k_1}\|_{Z_{k_1}})\cdot(2^{dk_2/2}\|f_{k_2}\|_{Z_{k_2}}),
\end{equation}
where $\mathcal{F}^{-1}(\widetilde{f}_{k_1})\in \{\mathcal{F}^{-1}(f_{k_1}),\overline{\mathcal{F}^{-1}(f_{k_1})}\}$.
\end{Lemmap1}

\begin{proof}[Proof of Lemma \ref{Lemmap1}] We may assume $k\leq k_2+20$. The bound \eqref{bt1} follows easily from \eqref{bt0} if $k_2\leq 99$ (compare with Case 1 below). Assume $k_2\geq 100$. In view of \eqref{mi2}, we may assume that
\begin{equation}\label{bt2}
f_{k_2}\text{ is supported in }I_{k_2}^{(d)}\times\R\cap \{(\xi_2,\tau_2):|\xi_2-v|\leq 2^{k_2-50}\}\text{ for some }v\in I_{k_2}^{(d)}.
\end{equation}
With $v$ as above, let $\widehat{v}=v/|v|\in\mathbb{S}^{d-1}$ and
\begin{equation*}
\widetilde{K}=\max(k_1+k_2,0)+100.
\end{equation*}
For $\e\in\{\e_1,\ldots,\e_L\}$ let
\begin{equation}\label{gg1}
\eta_{k,\e}(d)(\xi)=\begin{cases}
\eta_k^{(d)}(\xi)\cdot\eta^+_{[k-10,k+10]}(\xi\cdot\e)&\text{ if }k\geq 100;\\
\eta_k^{(d)}(\xi)&\text{ if }k< 100.
\end{cases}
\end{equation}
In view of  \eqref{mi2}, for \eqref{bt1} it suffices to prove that for any $\e\in\{\e_1,\ldots,\e_L\}$
\begin{equation}\label{bt1}
2^{dk/2}\|\eta_{k,\e}^{(d)}(\xi)\cdot (\widetilde{f}_{k_1}\ast f_{k_2})\|_{Z_k}\leq C2^{-|k_2-k|/4}(2^{dk_1/2}\|f_{k_1}\|_{Z_{k_1}})\cdot(2^{dk_2/2}\|f_{k_2}\|_{Z_{k_2}}),
\end{equation}

Using \eqref{qq1} with $j_1=j_2=0$, we estimate first
\begin{equation}\label{bt3}
\begin{split}
2^{dk/2}&\|\eta_{\leq\widetilde{K}-1}(\tau+|\xi|^2)\cdot \eta_{k,\e}^{(d)}(\xi)\cdot (\widetilde{f}_{k_1}\ast f_{k_2})\|_{Z_k}\\
&\leq C2^{dk/2}2^{\widetilde{K}/2}\beta_{k,\widetilde{K}}\|\eta_{k,\e}^{(d)}(\xi)\cdot (\widetilde{f}_{k_1}\ast f_{k_2})\|_{L^2}\\
&\leq C2^{dk/2}2^{\widetilde{K}/2}(1+2^{(\widetilde{K}-2k_+)/2})\cdot 2^{-\widetilde{K}/2}\cdot (2^{dk_1/2}\|f_{k_1}\|_{Z_{k_1}})\cdot(\|f_{k_2}\|_{Z_{k_2}})\\
&\leq C2^{d(k-k_2)/2}\cdot 2^{k_2-k}\cdot (2^{dk_1/2}\|f_{k_1}\|_{Z_{k_1}})\cdot(2^{dk_2/2}\|f_{k_2}\|_{Z_{k_2}}).
\end{split}
\end{equation}
It remains to estimate
\begin{equation}\label{bt4}
\begin{split}
2^{dk/2}\|\eta_{\geq \widetilde{K}}&(\tau+|\xi|^2)\cdot \eta_{k,\e}^{(d)}(\xi)\cdot (\widetilde{f}_{k_1}\ast f_{k_2})\|_{Z_k}\\
&\leq C2^{-|k_2-k|/4}(2^{dk_1/2}\|f_{k_1}\|_{Z_{k_1}})\cdot (2^{dk_2/2}\|f_{k_2}\|_{Z_{k_2}}).
\end{split}
\end{equation}
Using the atomic decomposition \eqref{mi1}, we analyze several cases.

{\bf{Case 1:}} $f_{k_2}=g_{k_2,j_2}\in X_{k_2}$, $f_{k_1}=g_{k_1,j_1}\in X_{k_1}$. We have to prove that
\begin{equation}\label{bt5}
\begin{split}
&2^{dk/2}\|\eta_{\geq\widetilde{K}}(\tau+|\xi|^2)\cdot \eta_{k,\e}^{(d)}(\xi)\cdot (\widetilde{g}_{k_1,j_1}\ast g_{k_2,j_2})\|_{Z_k}\\
&\leq C2^{-|k_2-k|/4}(2^{dk_1/2}2^{j_1/2}\beta_{k_1,j_1}\|g_{k_1,j_1}\|_{L^2})\cdot(2^{dk_2/2}2^{j_2/2}\beta_{k_2,j_2}\|g_{k_2,j_2}\|_{L^2}).
\end{split}
\end{equation}
We may assume that $\max(j_1,j_2)\geq \widetilde{K}-10$. Using \eqref{bt0}, the left-hand side of \eqref{bt5} is dominated by
\begin{equation*}
\begin{split}
C2^{dk/2}&2^{\jma/2}\beta_{k,\jma}\sup_{j\leq \jma+C}\|\mathbf{1}_{D_{k,j}}\cdot (\widetilde{g}_{k_1,j_1}\ast g_{k_2,j_2})\|_{L^2}\\
&\leq C2^{dk/2}2^{\jma/2}\beta_{k,\jma}\cdot 2^{dk_1/2}2^{\jmi/2}\|g_{k_1,j_1}\|_{L^2}\cdot \|g_{k_2,j_2}\|_{L^2},
\end{split}
\end{equation*}
which gives \eqref{bt5} using the simple inequality (see the definition \eqref{v1.1})
\begin{equation}\label{vi1.6}
\beta_{k,\jma}\leq C2^{k_2-k}\beta_{k_2,\jma}\leq C2^{k_2-k}\beta_{k_1,j_1}\beta_{k_2,j_2}.
\end{equation}

{\bf{Case 2:}} $f_{k_2}=g_{k_2,j_2}\in X_{k_2}$, $f_{k_1}\in Y^{\e_l}_{k_1}$, $l\in\{1,\ldots,L\}$. We have to prove that
\begin{equation}\label{bt6}
\begin{split}
2^{dk/2}\|&\eta_{\geq\widetilde{K}}(\tau+|\xi|^2)\cdot \eta_{k,\e}^{(d)}(\xi)\cdot (\widetilde{f}_{k_1}\ast g_{k_2,j_2})\|_{Z_k}\\
&\leq C2^{-|k_2-k|/4}(2^{dk_1/2}\|f_{k_1}\|_{Y^{\e_l}_{k_1}})\cdot(2^{dk_2/2}2^{j_2/2}\beta_{k_2,j_2}\|g_{k_2,j_2}\|_{L^2}).
\end{split}
\end{equation}
In view of the definitions, we may assume that $j_2\geq\widetilde{K}-10$ (otherwise the left-hand side of \eqref{bt6} is equal to $0$). We estimate the left-hand side of  \eqref{bt6} by
\begin{equation*}
\begin{split}
&C2^{dk/2}2^{j_2/2}\beta_{k,j_2}\|\widetilde{f}_{k_1}\ast g_{k_2,j_2}\|_{L^2}\leq C2^{dk/2}2^{j_2/2}\beta_{k,j_2}\cdot \|\mathcal{F}^{-1}(\widetilde{f}_{k_1})\|_{L^\infty}\cdot \|g_{k_2,j_2}\|_{L^2},
\end{split}
\end{equation*}
which gives \eqref{bt6}, in view of \eqref{lb44}.

{\bf{Case 3:}} $f_{k_2}=f_{k_2}^{\e_l,k'}\in Y_{k_2}^{\e_l,k'}$, $k'\leq k_1+20$, $f_{k_1}\in Y_{k_1}^{\e_{l'}}$, $l,l'\in\{1,\ldots, L\}$. In view of \eqref{pp2} and the analysis in Case 1 and Case 2 above, we may assume that $\widetilde{f}_{k_1}$ is supported in the set $\{(\xi_1,\tau_1):|\tau_1+|\xi_1|^2|\leq 2^{2k_1+50}\}$ and $f_{k_2}^{\e_l,k'}$ is supported in the set $\{(\xi_2,\tau_2):|\tau_2+|\xi_2|^2|\leq 2^{k+k'}\leq 2^{k+k_1+20}\}$. In this case the  left-hand side of  \eqref{bt4} is equal to $0$.

The analysis in Case 1, Case 2, and Case 3 suffices to prove \eqref{bt4} if $k_1\geq k_2-10$. So we may assume from now on that
\begin{equation}\label{as1}
k_1\leq k_2-10\text{ and }|k-k_2|\leq 2.
\end{equation}

{\bf{Case 4:}} $f_{k_2}=f_{k_2}^{\e_l,k'}\in Y_{k_2}^{\e_l,k'}$, $k'\in T_k$, $k'\leq k_1+20$, $f_{k_1}=g_{k_1,j_1}\in X_{k_1}$, $l\in\{1,\ldots, L\}$. In view of \eqref{pp2} and the analysis in Case 1 above, we may assume that $k'\geq 100$ and $f_{k_2}^{\e_l,k'}$ is supported in the set $\{(\xi_2,\tau_2):|\tau_2+|\xi_2|^2|\leq 2^{k+k'-100}\}$. Then we may assume $j_1\geq \widetilde{K}-10$; for \eqref{bt4} it suffices to prove that
\begin{equation}\label{bt10}
\begin{split}
2^{dk/2}2^{j_1/2}\beta_{k,j_1}&||\widetilde{g}_{k_1,j_1}\ast f_{k_2}^{\e_l,k'}||_{L^2}\\
&\leq C(2^{dk_1/2}2^{j_1/2}\beta_{k_1,j_1}\|g_{k_1,j_1}\|_{L^2})\cdot(2^{dk_2/2}\|f_{k_2}^{\e_l,k'}\|_{Y_{k_2}^{\e_l,k'}}).
\end{split}
\end{equation}
We use Lemma \ref{Lemmar1}, so we may assume
\begin{equation}\label{bt11}
\begin{cases}
&f^{\e_l,k'}_{k_2}(\xi_1\e_l+\xi',\tau)=2^{-k'/2}\cdot\frac{\eta_{\leq k'-100}(\xi_1-M)}{\xi_1-M+i/2^{k'}}\cdot h(\xi',\tau)\\
&||h||_{L^2_{\xi',\tau}}\leq C||f^{\e',k'}_{k_2}||_{Y_{k_2}^{\e_l,k'}},
\end{cases}
\end{equation}
where $\xi_1\in\R$, $\xi'\in P_{\e_l}$, $M=M(\xi',\tau)=(-\tau-|\xi'|^2)^{1/2}$, and $h$ is supported in
\begin{equation}\label{bt12}
S^{\e_l}_{k_2,k'}=\{(\xi',\tau)\in P_{\e_l}\times\R:-\tau-|\xi'|^2\in[2^{2k'-80},2^{2k'+10}],\,|\xi'|\leq  2^{k_2+1}\}.
\end{equation}
In view of \eqref{as1}, for \eqref{bt10}, it suffices to prove that for $h\in L^2(P_{\e_l}\times\R)$ and $f^{\e_l,k'}_{k_2}$ as in \eqref{bt11},
\begin{equation}\label{bt15}
||\widetilde{g}_{k_1,j_1}\ast f_{k_2}^{\e_l,k'}||_{L^2}\leq C2^{dk_1/2}\|g_{k_1,j_1}\|_{L^2}\cdot||h||_{L^2_{\xi',\tau}}.
\end{equation}
For later use, we prove \eqref{bt15} without using the restriction $k'\leq k_1+20$. We estimate the $L^2$ norm in the left-hand side of \eqref{bt15} by duality: the left-hand side of \eqref{bt15} is bounded by
\begin{equation*}
\begin{split}
2^{-k'/2}&\sup_{ \|a\|_{L^2}=1}\Big|\int_{(\R\times P_{\e_l}\times\R)^2}\widetilde{g}_{k_1,j_1}(\eta_1\e_l+\eta',\beta)\cdot h(\xi',\tau)\\
&\frac{\eta_{\leq k'-100}(\xi_1-M)}{\xi_1-M+i/2^{k'}}\cdot a(\xi_1+\eta_1,\xi'+\eta',\tau+\beta)\,d\xi_1d\eta_1d\xi'd\eta'd\tau d\beta\Big|.
\end{split}
\end{equation*}
Using the boundedness of the Hilbert transform on $L^2$ and then H\"{o}lder's ine\-qua\-li\-ty in the variables $(\xi',\tau,\beta)$, this is bounded by
\begin{equation*}
C2^{-k'/2}||h||_{L^2}\cdot 2^{k'/2}\cdot \int_{\R\times P_{\e_l}}\Big[\int_{\R}|\widetilde{g}_{k_1,j_1}(\eta_1\e_l+\eta',\beta)|^2\,d\beta\Big]^{1/2}\,d\eta_1d\eta',
\end{equation*}
which easily gives \eqref{bt15}.

{\bf{Case 5:}} $f_{k_2}=f_{k_2}^{\e_l,k'}\in Y_{k_2}^{\e_l,k'}$, $k'\geq k_1+20$, $f_{k_1}=f_{k_1,\leq\widetilde{K}}\in Z_{k_1}$, $l\in\{1,\ldots, L\}$. We may assume also $k'\geq  100$ and notice that
\begin{equation}\label{bt91}
f_{k_2}^{\e_l,k'}\ast \widetilde{f}_{k_1,\leq\widetilde{K}}\text{ is supported in the set }\{(\xi,\tau):\xi\cdot\e_l\in[2^{k'-2},2^{k'+2}]\}.
\end{equation}
We have to prove that
\begin{equation}\label{bt90}
\begin{split}
2^{dk/2}\|\eta_{\geq\widetilde{K}}&(\tau+|\xi|^2)\cdot \eta_{k,\e}^{(d)}(\xi)\cdot (\widetilde{f}_{k_1,\leq\widetilde{K}}\ast f_{k_2}^{\e_l,k'})\|_{Z_k}\\
&\leq C(2^{dk_1/2}\|f_{k_1,\leq\widetilde{K}}\|_{Z_{k_1}})\cdot (2^{dk_2/2}\|f_{k_2}^{\e_l,k'}\|_{Y^{\e_l,k'}_{k_2}}).
\end{split}
\end{equation}
We will use \eqref{bt70} implicitly in some of the estimates below, and write
\begin{equation}\label{bt94}
\begin{split}
-\mathcal{F}^{-1}&[(\tau+|\xi|^2+i)\cdot (\widetilde{f}_{k_1,\leq \widetilde{K}}\ast f_{k_2}^{\e_l,k'})]\\
&=(i\partial_t+\Delta_x-i)\mathcal{F}^{-1}(f_{k_2}^{\e_l,k'})\cdot \mathcal{F}^{-1}(\widetilde{f}_{k_1,\leq\widetilde{K}})\\
&+\mathcal{F}^{-1}(f_{k_2}^{\e_l,k'})\cdot (i\partial_t+\Delta_x)\mathcal{F}^{-1}(\widetilde{f}_{k_1,\leq \widetilde{K}})\\
&+2\nabla_x\mathcal{F}^{-1}(f_{k_2}^{\e_l,k'})\cdot \nabla_x\mathcal{F}^{-1}(\widetilde{f}_{k_1,\leq\widetilde{K}}).
\end{split}
\end{equation}
Thus, using \eqref{bt91},
\begin{equation}\label{qq20}
\begin{split}
2^{dk/2}&\|\eta_{\geq\widetilde{K}}(\tau+|\xi|^2)\cdot \eta_{k,\e}^{(d)}(\xi)\cdot (\widetilde{f}_{k_1,\leq\widetilde{K}}\ast f_{k_2}^{\e_l,k'})\|_{Z_k}\\
&\leq C2^{dk/2}2^{-k'/2}\gamma_{k_2,k'}\cdot ||(i\partial_t+\Delta_x-i)\mathcal{F}^{-1}(f_{k_2}^{\e_l,k'})\cdot \mathcal{F}^{-1}(\widetilde{f}_{k_1,\leq \widetilde{K}})||_{L^{1,2}_{\e_l}}\\
&+C 2^{dk/2}2^{-\widetilde{K}/2}||\mathcal{F}^{-1}(f_{k_2}^{\e_l,k'})\cdot (i\partial_t+\Delta_x)\mathcal{F}^{-1}(\widetilde{f}_{k_1,\leq \widetilde{K}})||_{L^2}\\\
&+C 2^{dk/2}2^{-\widetilde{K}/2}||\nabla_x\mathcal{F}^{-1}(f_{k_2}^{\e_l,k'})\cdot \nabla_x\mathcal{F}^{-1}(\widetilde{f}_{k_1,\leq \widetilde{K}})||_{L^2}.
\end{split}
\end{equation}
We estimate the first term in the right-hand side of \eqref{qq20} by
\begin{equation*}
C2^{dk/2}2^{-k'/2}\gamma_{k_2,k'}\cdot ||(i\partial_t+\Delta_x-i)\mathcal{F}^{-1}(f_{k_2}^{\e_l,k'})\|_{L^{1,2}_{\e_l}}\cdot \|\mathcal{F}^{-1}(\widetilde{f}_{k_1,\leq \widetilde{K}})||_{L^\infty},
\end{equation*}
which is bounded by the right-hand side of \eqref{bt90} in view of Lemma \ref{Lemmaa2}. We estimate the last two terms in the right-hand side of \eqref{qq20} by
\begin{equation*}
C 2^{dk/2}2^{-\widetilde{K}/2}\cdot 2^{\widetilde{K}}\|f_{k_2}^{\e_l,k'}\ast \widetilde{f}_{k_1,\leq \widetilde{K}}||_{L^2},
\end{equation*}
which is bounded by the right-hand side of \eqref{bt90}, using \eqref{qq1}.

{\bf{Case 6:}} $f_{k_2}=f_{k_2}^{\e_l,k'}\in Y_{k_2}^{\e_l,k'}$, $k'\geq k_1+20$, $f_{k_1}=g_{k_1,j_1}\in X_{k_1}$, $j_1\geq\widetilde{K}$, $l\in\{1,\ldots, L\}$. Then, using Lemma \ref{Lemmas1}, we decompose
\begin{equation*}
f_{k_2}^{\e_l,k'}=f_{k_2,\leq j_1-10}^{\e_l,k'}+f_{k_2,\geq j_1+10}^{\e_l,k'}+X_{k_2}.
\end{equation*}
In view of the analysis in Case 1, for \eqref{bt90} it suffices to prove that
\begin{equation}\label{bt95}
\begin{split}
&2^{dk/2}\|\eta_{k,\e}^{(d)}(\xi)\cdot (\widetilde{g}_{k_1,j_1}\ast f_{k_2,\leq j_1-10}^{\e_l,k'})\|_{Z_k}\\
&+2^{dk/2}\|\eta_{\geq j_1}(\tau+|\xi|^2)\cdot \eta_{k,\e}^{(d)}(\xi)\cdot (\widetilde{g}_{k_1,j_1}\ast f_{k_2,\geq  j_1+10}^{\e_l,k'})\|_{Z_k}\\
&\leq C(2^{dk_1/2}2^{j_1/2}\beta_{k_1,j_1}\|g_{k_1,j_1}\|_{L^2})\cdot (2^{dk_2/2}\|f_{k_2}^{\e_l,k'}\|_{Y^{\e_l,k'}_{k_2}}).
\end{split}
\end{equation}
The bound for the first term follows easily using the $L^2$ norm and  \eqref{bt15}.
To control the second term in the right-hand side of \eqref{bt95} we use again the decomposition \eqref{bt94}, as well as \eqref{bt91}, and estimate it (as in \eqref{qq20}) by
\begin{equation}\label{pp8}
\begin{split}
&C2^{dk/2}2^{-k'/2}\gamma_{k,k'}||(i\partial_t+\Delta_x-i)\mathcal{F}^{-1}(f_{k_2,\geq j_1+10}^{\e_l,k'})\cdot \mathcal{F}^{-1}(\widetilde{g}_{k_1,j_1})||_{L^{1,2}_{\e_l}}\\
&+C 2^{dk/2}2^{-j_1/2}||\mathcal{F}^{-1}(f_{k_2,\geq j_1+10}^{\e_l,k'})\cdot (i\partial_t+\Delta_x)\mathcal{F}^{-1}(\widetilde{g}_{k_1,j_1})||_{L^2}\\
&+C 2^{dk/2}2^{-j_1/2}||\nabla_x\mathcal{F}^{-1}(f_{k_2,\geq j_1+10}^{\e_l,k'})\cdot \nabla_x\mathcal{F}^{-1}(\widetilde{g}_{k_1,j_1})||_{L^2}.
\end{split}
\end{equation}
We estimate the first term in the right-hand side of \eqref{pp8} by
\begin{equation*}
C2^{dk/2}2^{-k'/2}\gamma_{k,k'}||(i\partial_t+\Delta_x-i)\mathcal{F}^{-1}(f_{k_2,\geq j_1+10}^{\e_l,k'})\|_{L^{1,2}_{\e_l}}\cdot \|\mathcal{F}^{-1}(\widetilde{g}_{k_1,j_1})||_{L^\infty},
\end{equation*}
which is bounded by the right-hand side of \eqref{bt95} in view of Lemma \ref{Lemmaa2}. We estimate the last two terms in the right-hand side of \eqref{pp8} by
\begin{equation*}
C 2^{dk/2}2^{-j_1/2}\cdot 2^{j_1}\|f_{k_2,\geq j_1+10}^{\e_l,k'}\ast \widetilde{g}_{k_1,j_1}||_{L^2},
\end{equation*}
which is bounded by the right-hand side of \eqref{bt90}, using \eqref{qq1}.
\end{proof}

\section{Dyadic bilinear estimates, II}\label{bilin2}

In this section we prove our second main bilinear estimate:

\newtheorem{Lemmap2}[Lemmap0]{Lemma}
\begin{Lemmap2}\label{Lemmap2}
If $k_1,k_2,k\in\Z$, $k_1\leq k_2-10$, $|k-k_2|\leq 2$, $f_{k_1}\in Z_{k_1}$, and $f_{k_2}\in Z_{k_2}$ then
\begin{equation}\label{bc1main}
\begin{split}
2^{dk/2}&\|\eta_k^{(d)}(\xi)\cdot (\tau+|\xi|^2+i)^{-1}\big[\widetilde{f}_{k_1}\ast [(\tau_2+|\xi_2|^2+i)f_{k_2}]\big]\|_{Z_k}\\
&\leq C(2^{dk_1/2}\|f_{k_1}\|_{Z_{k_1}})\cdot(2^{dk_2/2}\|f_{k_2}\|_{Z_{k_2}}),
\end{split}
\end{equation}
where $\mathcal{F}^{-1}(\widetilde{f}_{k_1})\in \{\mathcal{F}^{-1}(f_{k_1}),\overline{\mathcal{F}^{-1}(f_{k_1})}\}$.
\end{Lemmap2}

In view of \eqref{mi2}, we may assume that
\begin{equation}\label{km2}
f_{k_2}\text{ is supported in }I_{k_2}^{(d)}\times\R\cap \{(\xi_2,\tau_2):|\xi_2-v|\leq 2^{k_2-50}\}\text{ for some }v\in I_{k_2}^{(d)},
\end{equation}
and let $\widehat{v}=v/|v|$. With $\eta_{k,\e}$ defined as in \eqref{gg1},  $\e\in\{\e_1,\ldots,\e_L\}$, for \eqref{bc1main} it suffices to prove that
\begin{equation}\label{bc1}
\begin{split}
2^{dk/2}&\|\eta_{k,\e}^{(d)}(\xi)\cdot (\tau+|\xi|^2+i)^{-1}\big[\widetilde{f}_{k_1}\ast [(\tau_2+|\xi_2|^2+i)f_{k_2}]\big]\|_{Z_k}\\
&\leq C(2^{dk_1/2}\|f_{k_1}\|_{Z_{k_1}})\cdot(2^{dk_2/2}\|f_{k_2}\|_{Z_{k_2}}),
\end{split}
\end{equation}
for any $\e\in\{\e_1,\ldots,\e_L\}$. We consider again several cases.

{\bf{Case 1:}} $k_1\leq 100$ and $f_{k_2}=g_{k_2,j_2}\in X_{k_2}$, $j_2\geq k_1+k_2+100$. In this case we may assume $f_{k_1}=g_{k_1,j_1}\in X_{k_1}$. For \eqref{bc1} it suffices to prove that
\begin{equation}\label{bc2}
\begin{split}
2^{dk/2}2^{j_2}&\|\eta_{k,\e}^{(d)}(\xi)\cdot (\tau+|\xi|^2+i)^{-1}(\widetilde{g}_{k_1,j_1}\ast g_{k_2,j_2})\|_{Z_k}\\
&\leq C(2^{dk_1/2}2^{j_1/2}\beta_{k_1,j_1}\|g_{k_1,j_1}\|_{L^2})\cdot(2^{dk_2/2}2^{j_2/2}\beta_{k_2,j_2}\|g_{k_2,j_2}\|_{L^2}),
\end{split}
\end{equation}

If $|j_1-j_2|\geq 500$ then the left-hand side of \eqref{bc2} is dominated by
\begin{equation}\label{bh20}
\begin{split}
&C2^{dk/2}2^{j_2}2^{-\max(j_1,j_2)}\cdot 2^{\jma/2}\beta_{k,\jma}\|\widetilde{g}_{k_1,j_1}\ast g_{k_2,j_2}\|_{L^2}\\
&\leq C2^{dk/2}2^{j_2}2^{-\max(j_1,j_2)/2}\beta_{k,\jma}\cdot 2^{dk_1/2}2^{\jmi/2}\|g_{k_1,j_1}\|_{L^2}\cdot \|g_{k_2,j_2}\|_{L^2},
\end{split}
\end{equation}
using \eqref{bt0}, which suffices for \eqref{bc2} in view of \eqref{vi1.6}.

If $|j_1-j_2|\leq 500$ then, using \eqref{bt0}, the left-hand side of \eqref{bc2} is dominated by
\begin{equation*}
\begin{split}
&C2^{dk/2}2^{j_2}\sum_{j\leq j_2+C}2^{-j/2}\beta_{k,j}\cdot 2^{j/2}2^{dk_1/2}(\|g_{k_1,j_1}\|_{L^2}\cdot \|g_{k_2,j_2}\|_{L^2})\\
&\leq C2^{dk_2/2}2^{dk_1/2}2^{j_2}\cdot (j_2+\beta_{k_2,j_2})\cdot (\|g_{k_1,j_1}\|_{L^2}\cdot \|g_{k_2,j_2}\|_{L^2}),
\end{split}
\end{equation*}
which suffices for \eqref{bc2} since $\beta_{k_1,j_1}\approx 2^{j_1/2}\approx 2^{j_2/2}$.

{\bf{Case 2:}} $k_1\leq 100$ and $f_{k_2}=g_{k_2,j_2}\in X_{k_2}$, $j_2\leq k_1+k_2+100$. In this case we may assume $f_{k_1}=g_{k_1,j_1}\in X_{k_1}$. For \eqref{bc1} it suffices to prove that
\begin{equation}\label{bc9}
\begin{split}
2^{dk/2}2^{j_2}&\|\eta_{k,\e}^{(d)}(\xi)\cdot (\tau+|\xi|^2+i)^{-1}(\widetilde{g}_{k_1,j_1}\ast g_{k_2,j_2})\|_{Z_k}\\
&\leq C(2^{dk_1/2}2^{j_1/2}\beta_{k_1,j_1}\|g_{k_1,j_1}\|_{Z_{k_1}})\cdot(2^{dk_2/2}2^{j_2/2}\|g_{k_2,j_2}\|_{L^2}),
\end{split}
\end{equation}
In view of  \eqref{km2}, we may assume that $g_{k_2,j_2}$ is supported in $I_{k_2}^{(d)}\cap \{(\xi_2,\tau_2):|\xi_2-v|\leq 2^{k_2-40}\}$ for some $v\in I_{k_2}^{(d)}$. Then we estimate the left-hand side of \eqref{bc9} (using the $L^{1,2}_{\e}$ norm and \eqref{xx1}) by
\begin{equation*}
2^{dk/2}2^{j_2}\cdot 2^{-k/2}||\mathcal{F}^{-1}(\widetilde{g}_{k_1,j_1})||_{L^{2,\infty}_{\e}}\cdot ||\mathcal{F}^{-1}(g_{k_2,j_2})||_{L^{2,2}_{\e}},
\end{equation*}
which gives \eqref{bc9} in view of  Lemma \ref{Lemmaa1}.

{\bf{Case 3:}} $k_1\leq 100$ and $f_{k_2}=f_{k_2}^{\e_l,k'}\in Y^{\e_l,k'}_{k_2}$, $k'\in T_{k_2}$, $l\in\{1,\ldots,L\}$. In view of the analysis in Cases 1 and 2, and  \eqref{pp2}, we  may assume $k'\geq 200$. Then, using Lemma  \ref{Lemmaa2}, \eqref{xx1}, and the fact that $\widetilde{f}_{k_1}\ast [(\tau_2+|\xi_2|^2+i)f^{\e_l,k'}_{k_2}]$ is supported in the set $\{(\xi,\tau):\xi\cdot\e_l\in[2^{k'-2},2^{k'+2}]\}$, we estimate
\begin{equation*}
\begin{split}
&2^{dk/2}\|\eta_{k,\e}^{(d)}(\xi)\cdot (\tau+|\xi|^2+i)^{-1}\big[\widetilde{f}_{k_1}\ast [(\tau_2+|\xi_2|^2+i)f^{\e_l,k'}_{k_2}]\big]\|_{Z_k}\\
&\leq C2^{dk/2}2^{-k'/2}\gamma_{k,k'}||\mathcal{F}^{-1}(\widetilde{f}_{k_1})||_{L^\infty}\cdot \|\mathcal{F}^{-1}[(\tau_2+|\xi_2|^2+i)f^{\e_l,k'}_{k_2}]\|_{L^{1,2}_{\e_l}},
\end{split}
\end{equation*}
which suffices for \eqref{bc1}.

Thus, from now on we may assume
\begin{equation}\label{km1}
k_1\geq 100\text{ which implies that }k_2\geq 100.
\end{equation}
 As in the proof of Lemma \ref{Lemmap1}, let
\begin{equation*}
\widetilde{K}=k_1+k_2+100,
\end{equation*}
and define $f_{k_2,\leq\widetilde{K}-1}$ as  in \eqref{qq9}. Then, using Lemma \ref{Lemmaa1} and \eqref{xx1},
\begin{equation*}
\begin{split}
2^{dk/2}&\|\eta_{k,\e}^{(d)}(\xi)\cdot (\tau+|\xi|^2+i)^{-1}\big[\widetilde{f}_{k_1}\ast [(\tau_2+|\xi_2|^2+i)f_{k_2,\leq\widetilde{K}-1}]\big]\|_{Z_k}\\
&\leq C2^{dk/2}\cdot 2^{-k/2}||\mathcal{F}^{-1}(\widetilde{f}_{k_1})\cdot \mathcal{F}^{-1}[(\tau_2+|\xi_2|^2+i)f_{k_2,\leq\widetilde{K}-1}]||_{L^{1,2}_{\e}}\\
&\leq C2^{dk/2}\cdot 2^{-k/2}||\mathcal{F}^{-1}(\widetilde{f}_{k_1})||_{L^{2,\infty}_{\e}}\cdot ||\mathcal{F}^{-1}[(\tau_2+|\xi_2|^2+i)f_{k_2,\leq\widetilde{K}-1}]||_{L^{2,2}_{\e}}\\
&\leq C2^{dk/2}2^{-k/2}2^{(d-1)k_1/2}||f_{k_1}||_{Z_{k_1}}\cdot 2^{\widetilde{K}/2}||f_{k_2,\leq\widetilde{K}-1}||_{Z_{k_2}}.
\end{split}
\end{equation*}
Thus, for  \eqref{bc1} it suffices to  prove that
\begin{equation}\label{bc10}
\begin{split}
2^{dk/2}&\|\eta_{k,\e}^{(d)}(\xi)\cdot (\tau+|\xi|^2+i)^{-1}\big[\widetilde{f}_{k_1}\ast [(\tau_2+|\xi_2|^2+i)f_{k_2,\geq\widetilde{K}}]\big]\|_{Z_k}\\
&\leq C(2^{dk_1/2}\|f_{k_1}\|_{Z_{k_1}})\cdot(2^{dk_2/2}\|f_{k_2}\|_{Z_{k_2}}).
\end{split}
\end{equation}
To prove \eqref{bc10} we analyze several more cases:

{\bf{Case 4:}} $f_{k_2}=g_{k_2,j_2}\in X_{k_2}$, $f_{k_1}=g_{k_1,j_1}\in X_{k_1}$. We  may assume $j_2\geq \widetilde{K}$ and $g_{k_2,j_2}$ is supported in $I_{k_2}^{(d)}\cap \{(\xi_2,\tau_2):|\xi_2-v|\leq 2^{k_2-40}\}$ for some $v\in I_{k_2}^{(d)}$. If $|j_1-j_2|\geq 10$ then the same $L^2$ estimate as  in Case 1 (see \eqref{bh20}) gives the desired estimate. If $|j_1-j_2|\leq 10$, then we estimate the left-hand side of  \eqref{bc10} (using the $L^{1,2}_{\e}$ norm and \eqref{xx1}) by
\begin{equation*}
\begin{split}
&C2^{dk/2}2^{j_2}\cdot [(j_2-2k_2)_++1]2^{-k/2}||\mathcal{F}^{-1}(\widetilde{g}_{k_1,j_1})||_{L^{2,\infty}_{\e}}\cdot ||\mathcal{F}^{-1}(g_{k_2,j_2})||_{L^{2,2}_{\e}}\\
&\leq C\frac{2^{(j_2-k-k_1)/2}\cdot [(j_2-2k_2)_++1]}{\beta_{k_1,j_1}\cdot\beta_{k_2,j_2}}\cdot (2^{dk_1/2}\|f_{k_1}\|_{Z_{k_1}})\cdot(2^{dk_2/2}\|f_{k_2}\|_{Z_{k_2}}).
\end{split}
\end{equation*}
which  is controlled by the right-hand side of \eqref{bc10}.

{\bf{Case 5:}} $f_{k_2}=g_{k_2,j_2}\in X_{k_2}$, $f_{k_1}\in Y_{k_1}^{\e_l}$, $l\in\{1,\ldots,L\}$. We  may assume $j_2\geq \widetilde{K}$. Using the $L^2$ norm, we estimate the left-hand side of  \eqref{bc10} by
\begin{equation*}
C2^{dk/2}2^{j_2}\cdot 2^{-j_2}2^{j_2/2}\beta_{k,j_2}||\mathcal{F}^{-1}(\widetilde{f}_{k_1})||_{L^{\infty}}\cdot ||\mathcal{F}^{-1}(g_{k_2,j_2})||_{L^2},
\end{equation*}
which suffices, in view of Lemma  \ref{Lemmaa2}.

{\bf{Case 6:}} $f_{k_2}=f_{k_2}^{\e_l,k'}\in Y_{k_2}^{\e_l,k'}$, $k'\in T_{k_2}$, $f_{k_1}\in Y_{k_1}^{\e_{l'}}$, $l,l'\in\{1,\ldots,L\}$. In view of \eqref{pp2} and the analysis in Cases 4 and 5, we may assume $\widetilde{K}\leq k_2+k'$, so $k'\geq k_1+10$. Thus $f_{k_2}^{\e_l,k'}\ast \widetilde{f}_{k_1}$ is supported in the set $\{(\xi,\tau):\xi\cdot\e_l\in[2^{k'-2},2^{k'+2}]\}$, and we can estimate the left-hand side of \eqref{bc10} by
\begin{equation*}
\begin{split}
&C2^{dk/2}2^{-k'/2}\gamma_{k,k'}||\mathcal{F}^{-1}\big[\widetilde{f}_{k_1}\ast [(\tau_2+|\xi_2|^2+i)f^{\e_l,k'}_{k_2,\geq\widetilde{K}}]\big]||_{L^{1,2}_{\e_l}}\\
&\leq C2^{dk/2}2^{-k'/2}\gamma_{k,k'}||\mathcal{F}^{-1}(\widetilde{f}_{k_1})||_{L^\infty}\cdot ||\mathcal{F}^{-1}[(\tau_2+|\xi_2|^2+i)f^{\e_l,k'}_{k_2,\geq\widetilde{K}}]||_{L^{1,2}_{\e_l}},
\end{split}
\end{equation*}
which suffices, in view of  Lemma \ref{Lemmaa2}.

\section{Dyadic bilinear estimates, III}\label{bilin3}

In this section we prove our last dyadic bilinear estimate:

\newtheorem{Lemmap3}{Lemma}[section]
\begin{Lemmap3}\label{Lemmap3}
If $k_1,k_2,k\in\Z$, $k_1\leq k_2+10$, $f_{k_1}\in Z_{k_1}$, and $f_{k_2}\in Z_{k_2}$, then
\begin{equation}\label{br1main}
\begin{split}
2^{dk/2}&\|\eta_k^{(d)}(\xi)\cdot (\tau+|\xi|^2+i)^{-1}\big[ [(\tau_1+|\xi_1|^2+i)f_{k_1}]\ast f_{k_2}\big]\|_{Z_k}\\
&\leq C2^{-|k_2-k|/4}\cdot (2^{dk_1/2}\|f_{k_1}\|_{Z_{k_1}})\cdot(2^{dk_2/2}\|f_{k_2}\|_{Z_{k_2}}).
\end{split}
\end{equation}
\end{Lemmap3}

In view of \eqref{mi2}, we may assume that
\begin{equation}\label{nn1}
f_{k_2}\text{ is supported in }I_{k_2}^{(d)}\times\R\cap \{(\xi_2,\tau_2):|\xi_2-v|\leq 2^{k_2-50}\}\text{ for some }v\in I_{k_2}^{(d)},
\end{equation}
and let $\widehat{v}=v/|v|$. With $\eta_{k,\e}$ as  in \eqref{gg1}, for \eqref{br1main} it suffices  to prove that
\begin{equation}\label{br1}
\begin{split}
2^{dk/2}&\|\eta_{k,\e}^{(d)}(\xi)\cdot (\tau+|\xi|^2+i)^{-1}\big[ [(\tau_1+|\xi_1|^2+i)f_{k_1}]\ast f_{k_2}\big]\|_{Z_k}\\
&\leq C2^{-|k_2-k|/4}\cdot (2^{dk_1/2}\|f_{k_1}\|_{Z_{k_1}})\cdot(2^{dk_2/2}\|f_{k_2}\|_{Z_{k_2}}),
\end{split}
\end{equation}
for any  $\e\in\{\e_1,\ldots,\e_l\}$. We consider again several cases.

{\bf{Case 1:}} $k_1\geq 100$, $k_2\geq k_1+10d$, $f_{k_1}\in Z_{k_1}$, $f_{k_2}=g_{k_2,j_2}\in X_{k_2}$. We may assume $|k_2-k|\leq 2$ and let $g_{k_1,j_1}=f_{k_1}\cdot\eta_{j_1}(\tau+|\xi|^2)$, $j_1\in\Z_+$. Since $2^{j_1/2}\beta_{k_1,j_1}||g_{k_1,j_1}||_{L^2}\leq C||f_{k_1}||_{Z_{k_1}}$ (see Lemma \ref{Lemmas1}), for \eqref{br1} it suffices to prove that
\begin{equation}\label{nn2}
\begin{split}
&2^{dk/2}2^{j_1}\|\eta_{k,\e}^{(d)}(\xi)\cdot (\tau+|\xi|^2+i)^{-1}(g_{k_1,j_1}\ast g_{k_2,j_2})\|_{Z_k}\\
&\leq C(1+2^{k_1-j_1/2})^{-1}\cdot (2^{dk_1/2}2^{j_1/2}\beta_{k_1,j_1}\|g_{k_1,j_1}\|_{L^2})\cdot(2^{dk_2/2}2^{j_2/2}\beta_{k_2,j_2}\|g_{k_2,j_2}\|_{L^2}).
\end{split}
\end{equation}

We have several subcases depending on $j_1$ and $j_2$. Assume first that
\begin{equation}\label{nn4}
j_1\leq k_1+k_2+10\text{ and }j_2\leq k_1+k_2+20.
\end{equation}
For \eqref{nn2} it suffices to prove that (with $\e$ as in the function $\eta_{k,\e}^{(d)}$ in the left-hand side of  \eqref{nn2})
\begin{equation}\label{nn3}
\begin{split}
&2^{dk/2}2^{j_1}2^{-k/2}\|\mathcal{F}^{-1}(g_{k_1,j_1}\ast g_{k_2,j_2})\|_{L^{1,2}_{\e}}\\
&\leq C(1+2^{k_1-j_1/2})^{-1}\cdot (2^{dk_1/2}2^{j_1/2}\beta_{k_1,j_1}\|g_{k_1,j_1}\|_{L^2})\cdot(2^{dk_2/2}2^{j_2/2}\|g_{k_2,j_2}\|_{L^2}).
\end{split}
\end{equation}
We use the cutoff functions $\chi^{(1)}_{k_1,n_1}$ and $\chi^\e_{k_1,n'}$ defined in \eqref{nn7} and \eqref{nn6}) to decompose
\begin{equation}\label{nn9}
\begin{cases}
&g_{k_2,j_2}=\sum_{n_1\in\Xi^{(1)}_{k_1},\,n_1\in[2^{k_2-30},2^{k_2+10}]}\,\,\,\sum_{n'\in\Xi^\e_{k_1},\,|n'|\leq 2^{k_2+10}}g_{k_2,j_2}^{n_1,n'};\\
&g_{k_2,j_2}^{n_1,n'}(\xi_1\e+\xi',\tau)=g_{k_2,j_2}(\xi_1\e+\xi',\tau)\cdot\chi^{(1)}_{k_1,n_1}(\xi_1)\cdot\chi^\e_{k_1,n'}(\xi').\\
\end{cases}
\end{equation}
In view of \eqref{nn4}, we have the identity
\begin{equation}\label{nn10}
\begin{split}
&h^{n_1,n'}(\xi_1\e+\xi',\tau):=g_{k_1,j_1}\ast g_{k_2,j_2}^{n_1,n'}(\xi_1\e+\xi',\tau)\\
&=\widetilde{\chi}^{(1)}_{k_1,n_1}(\xi_1)\cdot\widetilde{\chi}^\e_{k_1,n'}(\xi')\cdot \eta_{\leq k_1+k_2+30d}(\tau+n_1^2+|n'|^2)\cdot h^{n_1,n'}(\xi_1\e+\xi',\tau)\\
&=\widetilde{\chi}^\e_{k_1,n'}(\xi')\cdot \eta_{\leq k_1+k_2+30d}(\tau+n_1^2+|n'|^2)\cdot h^{n_1,n'}(\xi_1\e+\xi',\tau).
\end{split}
\end{equation}
For simplicity of notation, let $\sum_{n_1,n'}$ denote the sum over $n_1$ and $n'$ as in \eqref{nn9}. Let
\begin{equation*}
\widetilde{h}^{n_1,n'}(x_1,\xi',\tau)=\int_{\R}e^{ix_1\xi_1}h^{n_1,n'}(\xi_1\e+\xi',\tau)\,d\xi_1,
\end{equation*}
thus, using \eqref{nn10},
\begin{equation*}
\int_{\R}e^{ix_1\xi_1}h^{n_1,n'}(\xi_1\e+\xi',\tau)\,d\xi_1=\widetilde{\chi}^\e_{k_1,n'}(\xi')\cdot \eta_{\leq k_1+k_2+30d}(\tau+n_1^2+|n'|^2)\cdot \widetilde{h}^{n_1,n'}(x_1,\xi',\tau).
\end{equation*}
We notice now that the supports in $(\xi',\tau)$ of $\widetilde{\chi}^\e_{k_1,n'}(\xi')\cdot \eta_{\leq k_1+k_2+30d}(\tau+n_1^2+|n'|^2)$ and $\widetilde{\chi}^\e_{k_1,m'}(\xi')\cdot \eta_{\leq k_1+k_2+30d}(\tau+m_1^2+|m'|^2)$ are disjoint unless $|n_1-m_1|+|n'-m'|\leq C2^{k_1}$ (recall that $n_1,m_1\approx 2^{k_2}$). Thus, for any $x_1\in\R$,
\begin{equation*}
\Big|\Big|\int_{\R}e^{ix_1\xi_1}\sum_{n_1,n'}h^{n_1,n'}\Big|\Big|^2_{L^2_{\xi',\tau}}\leq C\sum_{n_1,n'}\|\widetilde{h}^{n_1,n'}(x_1,\xi',\tau)\|_{L^2_{\xi',\tau}}^2.
\end{equation*}
Thus, using Plancherel theorem, the left-hand side of \eqref{nn3} is dominated by
\begin{equation}\label{nn20}
C2^{dk/2}2^{j_1}2^{-k/2}\int_{\R}\big[\sum_{n_1,n'}\|\mathcal{F}^{-1}(h^{n_1,n'})(x_1\e+x',t)\|_{L^2_{x',t}}^2\big]^{1/2}\,dx_1.
\end{equation}
We use now the definition of $h^{n_1,n'}$ in \eqref{nn10} to estimate, for any $x_1\in\R$,
\begin{equation*}
\begin{split}
\|\mathcal{F}^{-1}&(h^{n_1,n'})(x_1\e+x',t)\|_{L^2_{x',t}}\\
&\leq C\|\mathcal{F}^{-1}(g_{k_1,j_1})(x_1\e+x',t)\|_{L^2_{x',t}}\cdot \|\mathcal{F}^{-1}(g_{k_2,j_2}^{n_1,n'})(x_1\e+x',t)\|_{L^\infty_{x',t}}.
\end{split}
\end{equation*}
Thus, the expression in \eqref{nn20} is dominated by
\begin{equation*}
\begin{split}
C&2^{dk/2}2^{j_1}2^{-k/2}\\
&\int_{\R}\|\mathcal{F}^{-1}(g_{k_1,j_1})(x_1\e+x',t)\|_{L^2_{x',t}}\cdot \big[\sum_{n_1,n'}\|\mathcal{F}^{-1}(g_{k_2,j_2}^{n_1,n'})(x_1\e+x',t)\|_{L^\infty_{x',t}}^2\big]^{1/2}\,dx_1.
\end{split}
\end{equation*}
By H\"older's inequality in $x_1$, this is dominated by
\begin{equation}\label{nn60}
C2^{dk/2}2^{j_1}2^{-k/2}\|g_{k_1,j_1}\|_{L^2}\cdot \big[\sum_{n_1,n'}\|\mathcal{F}^{-1}(g_{k_2,j_2}^{n_1,n'})\|_{L^{2,\infty}_{\e}}^2\big]^{1/2}.
\end{equation}
Using the bound \eqref{mm1} in Lemma \ref{Lemmaa1}, this is dominated by
\begin{equation}\label{nn61}
C2^{dk_2/2}2^{j_1}2^{-k_2/2}\|g_{k_1,j_1}\|_{L^2}\cdot 2^{(d-1)k_1/2}\cdot 2^{(k_2-k_1)/2}\cdot (2^{j_2/2}\|g_{k_2,j_2}\|_{L^2}),
\end{equation}
which suffices for \eqref{nn3} since $\beta_{k_1,j_1}=1+2^{j_1/2-k_1}$.

Assume now that
\begin{equation}\label{nn30}
j_1\leq k_1+k_2+10\text{ and }j_2\geq  k_1+k_2+20.
\end{equation}
For \eqref{nn2} it suffices to prove that
\begin{equation}\label{nn31}
\begin{split}
&2^{dk/2}2^{j_1}2^{-j_2/2}\beta_{k,j_2}\|g_{k_1,j_1}\ast g_{k_2,j_2}\|_{L^2}\\
&\leq C(1+2^{k_1-j_1/2})^{-1}\cdot (2^{dk_1/2}2^{j_1/2}\beta_{k_1,j_1}\|g_{k_1,j_1}\|_{L^2})\cdot(2^{dk_2/2}2^{j_2/2}\beta_{k_2,j_2}\|g_{k_2,j_2}\|_{L^2}).
\end{split}
\end{equation}
Using Lemma \ref{Lemmap0}, we estimate the left-hand side of \eqref{nn31} by
\begin{equation*}
2^{dk_2/2}2^{j_1}2^{-j_2/2}\beta_{k,j_2}\cdot 2^{-j_2/2}(\beta_{k_1,j_1}\cdot\beta_{k_2,j_2})^{-1}\cdot (2^{dk_1/2}\|g_{k_1,j_1}\|_{Z_{k_1}})\cdot \|g_{k_2,j_2}\|_{Z_{k_2}},
\end{equation*}
which suffices for \eqref{nn31}. In this case we have proved the stronger bound
\begin{equation}\label{nn2.2}
\begin{split}
2^{dk/2}&2^{j_1}\|\eta_{k,\e}^{(d)}(\xi)\cdot (\tau+|\xi|^2+i)^{-1}(g_{k_1,j_1}\ast g_{k_2,j_2})\|_{Z_k}\leq C(1+2^{k_1-j_1/2})^{-1}\\
&\times 2^{(k_1-k_2)/2}\cdot (2^{dk_1/2}2^{j_1/2}\beta_{k_1,j_1}\|g_{k_1,j_1}\|_{L^2})\cdot(2^{dk_2/2}2^{j_2/2}\beta_{k_2,j_2}\|g_{k_2,j_2}\|_{L^2}).
\end{split}
\end{equation}

Assume now that
\begin{equation}\label{nn40}
j_1\geq k_1+k_2+10\text{ and }|j_2-j_1|\geq 10.
\end{equation}
Since the sequence $2^{-j/2}\beta_{k,j}$ is decreasing in $j$, for \eqref{nn2} it suffices to prove the stronger bound
\begin{equation}\label{nn41}
\begin{split}
&2^{dk/2}2^{j_1}2^{-j_1/2}\beta_{k,j_1}\cdot \sup_{j\in\Z_+}\|\mathbf{1}_{D_{k,j}}\cdot g_{k_1,j_1}\ast g_{k_2,j_2}\|_{L^2}\\
&\leq C2^{(k_1-k_2)/2}\cdot (2^{dk_1/2}2^{j_1/2}\beta_{k_1,j_1}\|g_{k_1,j_1}\|_{L^2})\cdot(2^{dk_2/2}2^{j_2/2}\beta_{k_2,j_2}\|g_{k_2,j_2}\|_{L^2}).
\end{split}
\end{equation}
Using \eqref{bt0}, we estimate the left-hand side of \eqref{nn41} by
\begin{equation*}
C2^{dk_2/2}2^{j_1/2}\beta_{k,j_1}\cdot 2^{dk_1/2}2^{j_2/2}\cdot\|g_{k_1,j_1}\|_{L^2}\cdot \|g_{k_2,j_2}\|_{L^2},
\end{equation*}
which suffices for \eqref{nn41}.

Finally, assume that
\begin{equation}\label{nn50}
j_1\geq k_1+k_2+10\text{ and }|j_2-j_1|\leq 10.
\end{equation}
Using \eqref{xx1}, for \eqref{nn2} it suffices to prove that
\begin{equation}\label{nn51}
\begin{split}
2^{dk/2}2^{j_1}&\cdot 2^{-k/2}[(j_2-2k_2)_++1]\|\mathcal{F}^{-1}(g_{k_1,j_1}\ast g_{k_2,j_2})\|_{L^{1,2}_\e}\\
&\leq C\cdot (2^{dk_1/2}2^{j_1/2}\beta_{k_1,j_1}\|g_{k_1,j_1}\|_{L^2})\cdot(2^{dk_2/2}2^{j_2/2}\beta_{k_2,j_2}\|g_{k_2,j_2}\|_{L^2}).
\end{split}
\end{equation}
Using \eqref{pr41}, we estimate the left-hand side of \eqref{nn51} by
\begin{equation*}
2^{dk_2/2}2^{j_1}2^{-k/2}[(j_2-2k_2)_++1]\cdot 2^{(d-1)k_1/2}2^{j_1/2}\|g_{k_1,j_1}\|_{L^2}\cdot \|g_{k_2,j_2}\|_{L^2},
\end{equation*}
which suffices for \eqref{nn51} since $\beta_{k_1,j_1}\geq 2^{(j_1-k_1-k_2)/2}$.

{\bf{Case 2:}} $k_1\geq 100$, $k_2\geq k_1+10d$, $f_{k_1}\in Z_{k_1}$, $f_{k_2}=f_{k_2}^{\e_l}\in Y^{\e_l}_{k_2}$, $l\in\{1,\ldots,L\}$. We may assume $|k_2-k|\leq 2$ and let $g_{k_1,j_1}=f_{k_1}\cdot\eta_{j_1}(\tau+|\xi|^2)$, $j_1\in\Z_+$. Since $2^{j_1/2}\beta_{k_1,j_1}||g_{k_1,j_1}||_{L^2}\leq C||f_{k_1}||_{Z_{k_1}}$, for \eqref{br1} it suffices to prove that
\begin{equation}\label{nm1}
\begin{split}
&2^{dk/2}2^{j_1}\|\eta_{k,\e}^{(d)}(\xi)\cdot (\tau+|\xi|^2+i)^{-1}(g_{k_1,j_1}\ast f_{k_2}^{\e_l})\|_{Z_k}\\
&\leq C(1+2^{k_1-j_1/2})^{-1}\cdot (2^{dk_1/2}2^{j_1/2}\beta_{k_1,j_1}\|g_{k_1,j_1}\|_{L^2})\cdot(2^{dk_2/2}\|f_{k_2}^{\e_l}\|_{Y_{k_2}^{\e_l}}).
\end{split}
\end{equation}

We consider two subcases. Assume first that
\begin{equation}\label{nm2}
j_1\leq k_1+k_2+10,
\end{equation}
and define $f_{k_2,\leq k_1+k_2+20}^{\e_l}$ and $f_{k_2,\geq k_1+k_2+21}^{\e_l}$ as  in \eqref{qq9}. To estimate
\begin{equation*}
2^{dk/2}2^{j_1}\|\eta_{k,\e}^{(d)}(\xi)\cdot (\tau+|\xi|^2+i)^{-1}(g_{k_1,j_1}\ast f_{k_2,\leq k_1+k_2+20}^{\e_l})\|_{Z_k}
\end{equation*}
we argue as in the proof of the bound  \eqref{nn4} in Case 1. The only difference is that  in passing from \eqref{nn60} to \eqref{nn61} we use the bound \eqref{mm1.1} in Lemma \ref{Lemmaa1}, instead of the bound \eqref{mm1}. To estimate
\begin{equation*}
2^{dk/2}2^{j_1}\|\eta_{k,\e}^{(d)}(\xi)\cdot (\tau+|\xi|^2+i)^{-1}(g_{k_1,j_1}\ast f_{k_2,\geq k_1+k_2+21}^{\e_l})\|_{Z_k}
\end{equation*}
we define $g_{k_2,j_2}=f_{k_2,\geq k_1+k_2+21}^{\e_l}\cdot \eta_{j_2}(\tau+|\xi|^2)$, $j_2\in[k_1+k_2+20,2k_2]\cap\Z$, and  use the bound \eqref{nn2.2} and Lemma \ref{Lemmas1}.

Assume now that
\begin{equation}\label{nm3}
j_1\geq k_1+k_2+10,
\end{equation}
and decompose
\begin{equation*}
f_{k_2}^{\e_l}=f^{\e_l}_{k_2,\leq j_1-10}+\sum_{j_2=j_1-9}^{2k_2}f^{\e_l}_{k_2}\cdot\eta_{j_2}(\tau+|\xi|^2).
\end{equation*}
The contribution  of the sum over $j_2$ in the expression above, which has  at  most $k_2-k_1$ terms, can be estimated using  \eqref{nn41} and  \eqref{nn51}. Then, we estimate the contribution of the function $f_{k_2,\leq j_1-10}^{\e_l}$ by
\begin{equation*}
C2^{dk/2}2^{j_1}\cdot 2^{-j_1/2}\beta_{k,j_1}||g_{k_1,j_1}\ast f_{k_2,\leq j_1-10}^{\e_l}||_{L^2}.
\end{equation*}
The bound \eqref{nm1} follows from \eqref{bt95}.

{\bf{Case 3:}} $k_2\leq C$. In this case, $k_1,k\leq C$ and we may assume $f_{k_1}=g_{k_1,j_1}\in X_{k_1}$ and $f_{k_2}=g_{k_2,j_2}\in X_{k_2}$. Since $\beta_{k_1,j_1}\approx 2^{j_1/2}$,  $\beta_{k_2,j_2}\approx 2^{j_2/2}$, $\beta_{k,j}\approx 2^{j/2}$, for  \eqref{br1} it suffices to prove that
\begin{equation}\label{kk1}
\begin{split}
2^{dk/2}2^{j_1}&\sum_{j\leq\max(j_1,j_2)+C}\|\mathbf{1}_{D_{k,j}}\cdot (g_{k_1,j_1}\ast g_{k_2,j_2})\|_{L^2}\\
&\leq C2^{-|k_2-k|/4}\cdot (2^{dk_1/2}2^{j_1}\|g_{k_1,j_1}\|_{L^2})\cdot(2^{dk_2/2}2^{j_2}\|g_{k_2,j_2}\|_{L^2}).
\end{split}
\end{equation}
This follows easily from \eqref{bt0}.

{\bf{Case 4:}} $k_1\leq 100$, $k_2\geq (k_1+10d)_+$, $f_{k_2}=g_{k_2,j_2}\in X_{k_2}$. We may assume $f_{k_1}=g_{k_1,j_1}\in X_{k_1}$, $\beta_{k_1,j_1}\approx 2^{j_1/2}$, and $|k_2-k|\leq 2$. For \eqref{br1} it suffices to prove that
\begin{equation}\label{bb1}
\begin{split}
&2^{dk/2}2^{j_1}\|\eta_{k,\e}^{(d)}(\xi)\cdot (\tau+|\xi|^2+i)^{-1}(g_{k_1,j_1}\ast g_{k_2,j_2})\|_{Z_k}\\
&\leq C(2^{dk_1/2}2^{j_1/2}2^{j_1/2}\|g_{k_1,j_1}\|_{L^2})\cdot(2^{dk_2/2}2^{j_2/2}\beta_{k_2,j_2}\|g_{k_2,j_2}\|_{L^2}).
\end{split}
\end{equation}
Assume first that
\begin{equation}\label{bb2}
j_1\leq k_1+k_2+10.
\end{equation}
Then, using \eqref{pr41} and \eqref{xx1}, we estimate the left-hand side of \eqref{bb1} by
\begin{equation}\label{cc2}
\begin{split}
C&2^{dk/2}2^{j_1}\cdot 2^{-k/2}[(j_2-2k_2)_++1]\|\mathcal{F}^{-1}(g_{k_1,j_1}\ast g_{k_2,j_2})\|_{L^{1,2}_\e}\\
&\leq C2^{dk_2/2}2^{j_1}2^{-k_2/2}[(j_2-2k_2)_++1]\cdot (2^{(d-1)k_1/2}2^{j_1/2}\|g_{k_1,j_1}\|_{L^2})\cdot \|g_{k_2,j_2}\|_{L^2},
\end{split}
\end{equation}
which suffices for \eqref{bb1}.

Assume now that
\begin{equation}\label{bb3}
j_1\geq  k_1+k_2+10\text{ and }|j_2-j_1|\geq 10.
\end{equation}
Since the sequence $2^{-j/2}\beta_{k,j}$ is decreasing in $j$, for \eqref{bb1} it suffices to prove that
\begin{equation}\label{bb4}
\begin{split}
&2^{dk/2}2^{j_1}2^{-j_1/2}\beta_{k,j_1}\cdot \sup_{j\in\Z_+}\|\mathbf{1}_{D_{k,j}}\cdot (g_{k_1,j_1}\ast g_{k_2,j_2})\|_{L^2}\\
&\leq C(2^{dk_1/2}2^{j_1/2}2^{j_1/2}\|g_{k_1,j_1}\|_{L^2})\cdot(2^{dk_2/2}2^{j_2/2}\beta_{k_2,j_2}\|g_{k_2,j_2}\|_{L^2}).
\end{split}
\end{equation}
Using \eqref{bt0}, we estimate the left-hand side of \eqref{bb4} by
\begin{equation*}
C2^{dk_2/2}2^{j_1/2}\beta_{k,j_1}\cdot 2^{dk_1/2}2^{j_2/2}\cdot\|g_{k_1,j_1}\|_{L^2}\cdot \|g_{k_2,j_2}\|_{L^2},
\end{equation*}
which suffices for \eqref{bb4}.

Finally, assume that
\begin{equation}\label{bb6}
j_1\geq k_1+k_2+10\text{ and }|j_2-j_1|\leq 10.
\end{equation}
For \eqref{bb1} it suffices to prove that
\begin{equation}\label{bb7}
\begin{split}
2^{dk/2}2^{j_1}&\sum_{j\leq j_1+20}2^{-j/2}\beta_{k,j}\cdot \|\mathbf{1}_{D_{k,j}}\cdot (g_{k_1,j_1}\ast g_{k_2,j_2})\|_{L^2}\\
&\leq C\cdot (2^{dk_1/2}2^{j_1/2}2^{j_1/2}\|g_{k_1,j_1}\|_{L^2})\cdot(2^{dk_2/2}2^{j_2/2}\beta_{k_2,j_2}\|g_{k_2,j_2}\|_{L^2}).
\end{split}
\end{equation}
Using \eqref{bt0}, we estimate the left-hand side of \eqref{nn51} by
\begin{equation*}
2^{dk_2/2}2^{j_1}\cdot \sum_{j\leq j_1+20}\beta_{k,j}\cdot 2^{dk_1/2}\|g_{k_1,j_1}\|_{L^2}\cdot \|g_{k_2,j_2}\|_{L^2},
\end{equation*}
which suffices for \eqref{bb7} since $\beta_{k,j}\leq C\beta_{k_2,j_2}$.

{\bf{Case 5:}} $k_1\leq 100$, $k_2\geq (k_1+10d)_+$, $f_{k_2}=f_{k_2}^{\e_l}\in Y^{\e_l}_{k_2}$. We may assume $f_{k_1}=g_{k_1,j_1}\in X_{k_1}$, $\beta_{k_1,j_1}\approx 2^{j_1/2}$, $k_2\geq 100$, and $|k_2-k|\leq 2$. For \eqref{br1} it suffices to prove that
\begin{equation}\label{cc1}
\begin{split}
&2^{dk/2}2^{j_1}\|\eta_{k,\e}^{(d)}(\xi)\cdot (\tau+|\xi|^2+i)^{-1}(g_{k_1,j_1}\ast f_{k_2}^{\e_l})\|_{Z_k}\\
&\leq C(2^{dk_1/2}2^{j_1/2}2^{j_1/2}\|g_{k_1,j_1}\|_{L^2})\cdot(2^{dk_2/2}\|f_{k_2}^{\e_l}\|_{Y_{k_2}^{\e_l}}).
\end{split}
\end{equation}
If $j_1\leq k_1+k_2+50$ then \eqref{cc1} follows using an estimate similar to \eqref{cc2} in Case 4 (clearly, $\|f_{k_2}^{\e_l}\|_{L^2}\leq \|f_{k_2}^{\e_l}\|_{Y_{k_2}^{\e_l}}$, using Lemma \ref{Lemmas1}). We assume
\begin{equation}\label{cc5}
j_1\geq k_1+k_2+50,
\end{equation}
and decompose
\begin{equation*}
f_{k_2}^{\e_l}=f_{k_2,\leq j_1-10}^{\e_l}+f_{k_2,\geq  j_1+10}^{\e_l}+X_{k_2}.
\end{equation*}
In view of \eqref{bt95},
\begin{equation*}
\begin{split}
&2^{dk/2}2^{j_1}\|\eta_{k,\e}^{(d)}(\xi)\cdot (\tau+|\xi|^2+i)^{-1}(g_{k_1,j_1}\ast f_{k_2,\leq j_1-10}^{\e_l})\|_{Z_k}\\
&\leq C(2^{dk_1/2}2^{j_1/2}2^{j_1/2}\|g_{k_1,j_1}\|_{L^2})\cdot(2^{dk_2/2}\|f_{k_2}^{\e_l}\|_{Y_{k_2}^{\e_l}}),
\end{split}
\end{equation*}
as desired. In addition,
\begin{equation*}
\begin{split}
&2^{dk/2}2^{j_1}\|\eta_{k,\e}^{(d)}(\xi)\cdot (\tau+|\xi|^2+i)^{-1}(g_{k_1,j_1}\ast f_{k_2,\geq j_1+10}^{\e_l})\|_{Z_k}\\
&\leq C2^{dk_2/2}2^{j_1}\sum_{j=j_1}^{2k+10}2^{-j/2}\|\mathbf{1}_{D_{k,j}}\cdot (g_{k_1,j_1}\ast f_{k_2,\geq  j_1+10}^{\e_l})\|_{L^2}\\
&\leq C2^{dk_2/2}2^{j_1}\big(\sum_{j=j_1}^{2k+10}2^{-j/2}\big)\cdot 2^{-j_1/2}\cdot 2^{-j_1/2}\cdot (2^{dk_1/2}\|g_{k_1,j_1}\|_{Z_{k_1}})\cdot( \|f_{k_2}^{\e_l}\|_{Y_{k_2}^{\e_l}}),
\end{split}
\end{equation*}
using Lemma \ref{Lemmap0}, since $j_1\leq 2k_2$. This completes the proof of \eqref{cc1}.

{\bf{Case 6:}} $k_1,k_2\geq 100d$, $|k_1-k_2|\leq 10d$, $f_{k_1}\in Z_{k_1}$, $f_{k_2}=g_{k_2,j_2}\in X_{k_2}$. Let $g_{k_1,j_1}=f_{k_1}\cdot \eta_{j_1}(\tau+|\xi|^2)$, $j_1\in\Z_+$. Since $2^{j_1/2}\beta_{k_1,j_1}\|g_{k_1,j_1}\|_{L^2}\leq C\|f_{k_1}\|_{Z_{k_1}}$, for \eqref{br1} it suffices to prove that
\begin{equation}\label{ll1}
\begin{split}
&2^{dk/2}2^{j_1}\|\eta_{k,\e}^{(d)}(\xi)\cdot (\tau+|\xi|^2+i)^{-1}(g_{k_1,j_1}\ast g_{k_2,j_2})\|_{Z_k}\leq C(1+2^{k_1-j_1/2})^{-1}\\
&\times 2^{-|k_2-k|/4}(2^{dk_1/2}2^{j_1/2}\beta_{k_1,j_1}\|g_{k_1,j_1}\|_{L^2})\cdot(2^{dk_2/2}2^{j_2/2}\beta_{k_2,j_2}\|g_{k_2,j_2}\|_{L^2}).
\end{split}
\end{equation}
Using the definition, we estimate the left-hand side of \eqref{ll1} by
\begin{equation}\label{ll2}
\begin{split}
&C2^{dk/2}2^{j_1}\|\eta_{k,\e}^{(d)}(\xi)\cdot \eta_{\leq 2k-201}(\tau+|\xi|^2)\cdot (\tau+|\xi|^2+i)^{-1}(g_{k_1,j_1}\ast g_{k_2,j_2})\|_{Y_k^\e}\\
&+C2^{dk/2}2^{j_1}\|\eta_{k,\e}^{(d)}(\xi)\cdot \eta_{\geq 2k-200}(\tau+|\xi|^2)\cdot (\tau+|\xi|^2+i)^{-1}(g_{k_1,j_1}\ast g_{k_2,j_2})\|_{Z_k}.
\end{split}
\end{equation}
We estimate the first term in \eqref{ll2} (which is nontrivial only if $k\geq 100$) by
\begin{equation*}
\begin{split}
&C2^{dk/2}2^{j_1}2^{-k/2}\|\mathcal{F}^{-1}(g_{k_1,j_1})\|_{L^2}\cdot \|\mathcal{F}^{-1}(g_{k_2,j_2})\|_{L^{2,\infty}_\e}\\
&\leq C2^{dk/2}2^{j_1}2^{-k/2}\|g_{k_1,j_1}\|_{L^2}\cdot 2^{(d-1)k_2/2}\|g_{k_2,j_2}\|_{Z_{k_2}},
\end{split}
\end{equation*}
using Lemma \ref{Lemmaa1}. This suffices for \eqref{ll1} since $\beta_{k_1,j_1}=1+2^{j_1/2-k_1}$. For the second term in \eqref{ll2} we use $L^2$ estimates. Assume first that
\begin{equation}\label{ll3}
j_1\leq 2k_1+30d.
\end{equation}
Then the second term in \eqref{ll2} is bounded by
\begin{equation*}
C2^{dk/2}2^{j_1}\sum_{j}2^{-j/2}\beta_{k,j}\|\mathbf{1}_{D_{k,j}}\cdot (g_{k_1,j_1}\ast g_{k_2,j_2})\|_{L^2},
\end{equation*}
where the sum is over $j\geq \max(0,2k-200)$ and $j\leq \max(2k_2,j_2)+C$. Since $\beta_{k,j}\approx 2^{j/2-k_+}$, using Lemma \ref{Lemmap0} this expression is bounded by
\begin{equation*}
\begin{split}
C2^{dk/2}2^{j_1}2^{-k_+}&[ |k_2-k_+|+(j_2-2k_2)_++1]\cdot 2^{-k_2}\cdot\beta_{k_2,j_2}^{-1}\\
&\times ( 2^{dk_1/2}\|g_{k_1,j_1}\|_{Z_{k_1}})\cdot \|g_{k_2,j_2}\|_{Z_{k_2}},
\end{split}
\end{equation*}
which suffices for \eqref{ll1} (recall that $d\geq 3$).

Assume now that
\begin{equation}\label{ll5}
j_1\geq 2k_1+30d\text{ and }|j_1-j_2|\leq 10.
\end{equation}
Then the second term in \eqref{ll2} is bounded by
\begin{equation*}
C2^{dk/2}2^{j_1}\sum_{j}2^{-j/2}\beta_{k,j}\|\mathbf{1}_{D_{k,j}}\cdot (g_{k_1,j_1}\ast g_{k_2,j_2})\|_{L^2},
\end{equation*}
where the sum is over $j\geq \max(0,2k-200)$ and $j\leq j_2+C$. Since $\beta_{k,j}\approx 2^{j/2-k_+}$, using Lemma \ref{Lemmap0} this expression is bounded by
\begin{equation*}
\begin{split}
C2^{dk/2}2^{j_1}2^{-k_+}&( |j_1-2k_+|+1)\cdot 2^{-j_2/2}\cdot(\beta_{k_1,j_1}\beta_{k_2,j_2})^{-1}\\
&\times ( 2^{dk_1/2}\|g_{k_1,j_1}\|_{Z_{k_1}})\cdot \|g_{k_2,j_2}\|_{Z_{k_2}},
\end{split}
\end{equation*}
which suffices for \eqref{ll1}.

Finally, assume that
\begin{equation}\label{ll10}
j_1\geq 2k_1+30d\text{ and }|j_1-j_2|\geq  10.
\end{equation}
Since $2^{-j/2}\beta_{k,j}\approx 2^{-k_+}$ for $j\geq 2k_+$, the second term in \eqref{ll2} is  bounded by
\begin{equation*}
\begin{split}
&C2^{dk/2}2^{j_1}2^{-k_+}\|g_{k_1,j_1}\ast g_{k_2,j_2}\|_{L^2}\\
&\leq C2^{dk/2}2^{j_1}2^{-k_+}\cdot 2^{-j_1/2}\beta_{k_1,j_1}^{-1}\cdot( 2^{dk_1/2}\|g_{k_1,j_1}\|_{Z_{k_1}})\cdot \|g_{k_2,j_2}\|_{Z_{k_2}},
\end{split}
\end{equation*}
using Lemma \ref{Lemmap0}, which suffices for \eqref{ll1}.

{\bf{Case 7:}} $k_1,k_2\geq 100d$, $|k_1-k_2|\leq 10d$, $f_{k_1}\in Z_{k_1}$, $f_{k_2}=f_{k_2}^{\e_l}\in Y^{\e_l}_{k_2}$, $l\in\{1,\ldots,L\}$. Let $g_{k_1,j_1}=f_{k_1}\cdot \eta_{j_1}(\tau+|\xi|^2)$, $j_1\in\Z_+$. Since $2^{j_1/2}\beta_{k_1,j_1}\|g_{k_1,j_1}\|_{L^2}\leq C\|f_{k_1}\|_{Z_{k_1}}$, for \eqref{br1} it suffices to prove that
\begin{equation}\label{dd1}
\begin{split}
&2^{dk/2}2^{j_1}\|\eta_{k,\e}^{(d)}(\xi)\cdot (\tau+|\xi|^2+i)^{-1}(g_{k_1,j_1}\ast f_{k_2}^{\e_l})\|_{Z_k}\leq C(1+2^{k_1-j_1/2})^{-1}\\
&\times 2^{-|k_2-k|/4}(2^{dk_1/2}2^{j_1/2}\beta_{k_1,j_1}\|g_{k_1,j_1}\|_{L^2})\cdot(2^{dk_2/2}\|f_{k_2}^{\e_l}\|_{Y_{k_2}^{\e_l}}).
\end{split}
\end{equation}
Using the definition, we estimate the left-hand side of \eqref{dd1} by
\begin{equation}\label{dd2}
\begin{split}
&C2^{dk/2}2^{j_1}\|\eta_{k,\e}^{(d)}(\xi)\cdot \eta_{\leq 2k-201}(\tau+|\xi|^2)\cdot (\tau+|\xi|^2+i)^{-1}(g_{k_1,j_1}\ast f_{k_2}^{\e_l})\|_{Y_k^\e}\\
&+C2^{dk/2}2^{j_1}\|\eta_{k,\e}^{(d)}(\xi)\cdot \eta_{\geq 2k-200}(\tau+|\xi|^2)\cdot (\tau+|\xi|^2+i)^{-1}(g_{k_1,j_1}\ast f_{k_2}^{\e_l})\|_{Z_k}.
\end{split}
\end{equation}
We estimate the first term in \eqref{dd2} (which is nontrivial only if $k\geq 100$) by
\begin{equation*}
\begin{split}
&C2^{dk/2}2^{j_1}2^{-k/2}\|\mathcal{F}^{-1}(g_{k_1,j_1})\|_{L^2}\cdot \|\mathcal{F}^{-1}(f_{k_2}^{\e_l})\|_{L^{2,\infty}_\e}\\
&\leq C2^{dk/2}2^{j_1}2^{-k/2}\|g_{k_1,j_1}\|_{L^2}\cdot 2^{(d-1)k_2/2}\|f_{k_2}^{\e_l}\|_{Y^{\e_l}_{k_2}},
\end{split}
\end{equation*}
using Lemma \ref{Lemmaa1}. This suffices for \eqref{dd1} since $\beta_{k_1,j_1}=1+2^{j_1/2-k_1}$. For the second term in \eqref{dd2} we use $L^2$ estimates. Assume first that
\begin{equation}\label{dd3}
j_1\leq 2k_1+30d.
\end{equation}
Then the second term in \eqref{dd2} is bounded by
\begin{equation*}
C2^{dk/2}2^{j_1}\sum_{j}2^{-j/2}\beta_{k,j}\|\mathbf{1}_{D_{k,j}}\cdot (g_{k_1,j_1}\ast f^{\e_l}_{k_2})\|_{L^2},
\end{equation*}
where the sum is over $j\geq \max(0,2k-200)$ and $j\leq 2k_2+C$. Since $\beta_{k,j}\approx 2^{j/2-k_+}$, using Lemma \ref{Lemmap0} this expression is bounded by
\begin{equation*}
\begin{split}
C2^{dk/2}2^{j_1}2^{-k_+}&[ |k_2-k_+|+1]\cdot 2^{-k_2}\cdot( 2^{dk_1/2}\|g_{k_1,j_1}\|_{Z_{k_1}})\cdot \|f_{k_2}^{\e_l}\|_{Z_{k_2}},
\end{split}
\end{equation*}
which suffices for \eqref{dd1} (recall that $d\geq 3$).

Assume now that
\begin{equation}\label{dd10}
j_1\geq 2k_1+30d.
\end{equation}
Since $2^{-j/2}\beta_{k,j}\approx 2^{-k_+}$ for $j\geq 2k_+$, the second term in \eqref{dd2} is  bounded by
\begin{equation*}
\begin{split}
&C2^{dk/2}2^{j_1}2^{-k_+}\|g_{k_1,j_1}\ast f_{k_2}^{\e_l}\|_{L^2}\\
&\leq C2^{dk/2}2^{j_1}2^{-k_+}\cdot 2^{-j_1/2}\beta_{k_1,j_1}^{-1}\cdot\|g_{k_1,j_1}\|_{Z_{k_1}}\cdot ( 2^{dk_2/2}\|g_{k_2,j_2}\|_{Z_{k_2}}),
\end{split}
\end{equation*}
using Lemma \ref{Lemmap0}, which suffices for \eqref{dd1}.

\section{A dyadic trilinear estimate}\label{trilin}

In this section we  prove the following trilinear estimate:

\newtheorem{Lemmaq1}{Lemma}[section]
\begin{Lemmaq1}\label{Lemmaq1}
If $k_1,k_2,k_3,k\in\Z$, $f_{k_1}\in Z_{k_1},\,f_{k_2}\in Z_{k_2},\,f_{k_3}\in Z_{k_3}$, and
\begin{equation}\label{vm0}
\min(k,k_2,k_3)\leq k_1+20,
\end{equation}
then
\begin{equation}\label{vm1}
\begin{split}
&2^{k_2+k_3}\cdot 2^{dk/2}\|\eta_{k}^{(d)}(\xi)\cdot (\tau+|\xi|^2+i)^{-1}\cdot(\widetilde{f}_{k_1}\ast\widetilde{f}_{k_2}\ast \widetilde{f}_{k_3})\|_{Z_k}\\
&\leq C2^{-|\max(k_1,k_2,k_3)-k|/4}\cdot (2^{dk_1/2}\|f_{k_1}\|_{Z_{k_1}})\cdot(2^{dk_2/2}\|f_{k_2}\|_{Z_{k_2}})\cdot (2^{dk_3/2}\|f_{k_3}\|_{Z_{k_3}}),
\end{split}
\end{equation}
where $\mathcal{F}^{-1}(\widetilde{f}_{k_l})\in \{\mathcal{F}^{-1}(f_{k_l}),\overline{\mathcal{F}^{-1}(f_{k_l})}\}$, $l=1,2,3$.
\end{Lemmaq1}

By symmetry, we may assume $k_2\leq k_3$. We start with the
following simple geometric observation: if
$\widehat{w}_1,\widehat{w}_2\in\mathbb{S}^{d-1}$ then there is
$\e\in\{\e_1,\ldots,\e_L\}$ such that
\begin{equation}\label{vm2}
\e\cdot \widehat{w}_1\geq 2^{-5}\text{  and }|\e\cdot
\widehat{w}_2|\geq 2^{-5}.
\end{equation}
To prove this, we may assume $\widehat{w}_1\cdot\widehat{w}_2\geq
0$  (by possibly replacing $\widehat{w}_2$ with $-\widehat{w}_2$)
and take $\e\in\{\e_1,\ldots,\e_L\}$ with
$\big|\e-(\widehat{w}_1+\widehat{w}_2)/|\widehat{w}_1+\widehat{w}_2|\,\big|\leq
2^{-100}$ (compare with  \eqref{vm4}).

The bound \eqref{xx1} shows that if $k\in\Z$ and
\begin{equation*}
f\text{ is supported in }I_k^{(d)}\times\R\cap
\{(\xi,\tau):\xi\cdot\e\geq 2^{k-40}\}\text{ for some }\e\in
\{\e_1,\ldots,\e_L\},
\end{equation*}
then, for any $J\geq 0$,
\begin{equation}\label{vm10}
||f\cdot\eta_{\leq J}(\tau+|\xi|^2)||_{Z_k}\leq
C\big((J-2k_+)_++1\big)\cdot 2^{-k/2}||\mathcal{F}^{-1}
[(\tau+|\xi|^2+i)\cdot f]||_{L^{1,2}_\e}.
\end{equation}

In view of \eqref{mi2}, we may assume that for $i=1,2,3$
\begin{equation}\label{vm20}
f_{k_i}\text{ is supported in }I_{k_i}^{(d)}\times\R\cap
\{(\xi,\tau):|\xi-v_i| \leq 2^{k_i-50}\}\text{ for some }v_i\in
I_{k_i}^{(d)},
\end{equation}
and it suffices to prove that for any $v\in I_{k}^{(d)}$
\begin{equation}\label{vm21}
\begin{split}
&2^{k_2+k_3}\cdot 2^{dk/2}\|\eta_{k}^{(d)}(\xi)\cdot \eta_0(|\xi-v|/2^{k-50})
\cdot (\tau+|\xi|^2+i)^{-1}\cdot(\widetilde{f}_{k_1}\ast\widetilde{f}_{k_2}
\ast \widetilde{f}_{k_3})\|_{Z_k}\\
&\leq C2^{-|\max(k_1,k_2,k_3)-k|/4}\cdot (2^{dk_1/2}\|f_{k_1}\|_{Z_{k_1}})
\cdot(2^{dk_2/2}\|f_{k_2}\|_{Z_{k_2}})\cdot (2^{dk_3/2}\|f_{k_3}\|_{Z_{k_3}}).
\end{split}
\end{equation}

Assume $k_1\leq k_3$ (in the case $k_1\geq k_3$ the bound
\eqref{vm41} below still holds, by a similar argument). Let
$\widehat{w}_1=v/|v|$, $\widehat{w}_2=v_3/|v_3|$, and fix $\e$ as
in \eqref{vm2}. Fix
\begin{equation}\label{vm40}
J=2\max(k_1,k_2,k_3,0)+100.
\end{equation}
Let $$F(\xi,\tau)=\eta_{k}^{(d)}(\xi)\cdot
\eta_0(|\xi-v|/2^{k-50}) \cdot
(\tau+|\xi|^2+i)^{-1}\cdot(\widetilde{f}_{k_1}\ast\widetilde{f}_{k_2}
\ast \widetilde{f}_{k_3})$$ and
$$\Pi=(2^{dk_1/2}\|f_{k_1}\|_{Z_{k_1}})
\cdot(2^{dk_2/2}\|f_{k_2}\|_{Z_{k_2}})\cdot
(2^{dk_3/2}\|f_{k_3}\|_{Z_{k_3}}).$$ Using \eqref{vm10}, Lemma
\ref{Lemmaa1} and Lemma \ref{Lemmas2},
\begin{equation}\label{vm41}
\begin{split}
&2^{k_2+k_3}\cdot 2^{dk/2}\|\eta_{\leq J}(\tau+|\xi|^2)\cdot
F\|_{Z_k}\leq C2^{k_2+k_3}2^{dk/2}\\
&\times
2^{-k/2}\big((J-2k_+)_++1\big)||\mathcal{F}^{-1}(\widetilde{f}_{k_1})||_
{L^{2,\infty}_\e}||\mathcal{F}^{-1}(\widetilde{f}_{k_2})||_
{L^{2,\infty}_\e}||\mathcal{F}^{-1}(\widetilde{f}_{k_3})||_
{L^{\infty,2}_\e}\\
&\leq C2^{k_2+k_3}2^{dk/2}2^{-k/2}\big((J-2k_+)_++1\big)\cdot
2^{-(k_1+k_2+k_3)/2}2^{-d\max(k_1,k_2,k_3)/2}\cdot \Pi,
\end{split}
\end{equation}
which is dominated by the right-hand side of \eqref{vm21},
provided that \eqref{vm0} holds and $d\geq 3$.

It remains to bound
\begin{equation}\label{vm50}
2^{k_2+k_3}\cdot 2^{dk/2}\|\eta_{\geq J+1}(\tau+|\xi|^2)\cdot
F\|_{Z_k}.
\end{equation}
We use the atomic decomposition \eqref{mi1} for the functions
$f_{k_1}$, $f_{k_2}$, and $f_{k_3}$, and notice that the
expression in \eqref{vm50} is not equal to $0$ only if at least
one of the functions $f_{k_1}$, $f_{k_2}$, or $f_{k_3}$ has
modulation $(\tau+|\xi|^2)\geq 2^{J-10}$. Let $J'\geq J-10$ denote
the highest of these modulations. Then we estimate
$2^{k_2+k_3}\cdot 2^{dk/2}\|\eta_{\leq J'+10}(\tau+|\xi|^2)\cdot
F\|_{Z_k}$ as in \eqref{vm41}. Using \eqref{pr41} or \eqref{gu41}
for the function with the high modulation, the right-hand side of
\eqref{vm41} is multiplied by at most
\begin{equation*}
C\beta_{\max(k_1,k_2,k_3),J'}^{-1}\cdot (|J'-J|+1)\leq C,
\end{equation*}
which suffices to complete the proof of \eqref{vm21}.

\end{document}